\input amstex
\documentstyle{amsppt}

\font\twolverm=cmssbx10 at 12pt
\font\twoolverm=cmssbx10 at 20pt
\font\twoplverm=cmssbx10 at 14pt

\font\bbol =cmssbx10 
\font\ita=cmssi10
\font\es=cmss10 
\font\esp=cmss10 at 6pt 
\input epsf
\define\te#1{{\text{\es#1}}}
\define\QQ{{\overline{\bold{ Q}}}}
\define\Fp{{\overline{\bold F}_p}}

%\magnification=1200

\hcorrection{1.5cm}

\nologo
\nopagenumbers
\topmatter\endtopmatter

\ 

\vskip1cm

\noindent{\twoolverm Galois action on diameter four trees}

\vskip1cm

\noindent{\twolverm Leonardo Zapponi}

\vskip.1cm

\noindent{\ita Max-Planck-Institut f\"ur Mathematik,}
{\ita Vivatsgasse 7, 53111 Bonn (Germany).\newline}
{\ita e-mail: zapponi\@mpim-bonn.mpg.de}

\vskip.4cm

\noindent{\bbol Abstract.}\es This paper is an extended and updated version of a talk presented by the author during the conference ``Th\'eorie de Galois et g\'eom\'etrie'', yeld in June 2001 at the C.I.R.M of Luminy. Its main object is the study of a class of dessins d'enfants, the so-called diameter four trees. These objects, first introduced by G. Shabat in [Sh], can be considered as the simplest non trivial example of \'etale covers of the projective line minus three points. Their arithmetic properties are still mysterious, and their study can inspire the understanding of  more general situations. Here, the main interest  is devoted to the action of the absolute Galois group on these (isomorphism classes of) coverings. In particular, in many cases,  we are able to distinguish Galois orbits and to describe the action of the decomposition groups. One other central result concerns  the study of wild ramification, for which we show how to reduce to the tame case, and then deduce some detailed arithmetical informations.

\baselineskip=14pt

\vskip1.4cm

\es 

\noindent{\twoplverm Contents}

\vskip.3cm 

1. Diameter four trees.

2. Standard and normalized models. 

3. Galois action, fields of moduli.

4. Fundamental equations for standard models.

5. $p$-congruence

6. Lifting normalized models from positive characteritic.

7. Good reduction, action of the decomposition groups.

8. Wild ramification above zero, Kummer models.

9. Wild ramification above infinity.

\vskip2cm

\specialhead\twoplverm 1 Diameter four trees\endspecialhead

\vskip.3cm

Let $K$ be an algebraically closed field and consider $n$ positive integers $a_1\leq\dots\leq a_n$. A {\bbol diameter four tree} over $K$, {\bbol of type} $(a_1,\dots,a_n)$, is an isomorphism class $\Cal T=[\beta]$ of coverings $\beta:\bold P^1_K@>\quad>>\bold P^1_K$ having the following ramification data:

\vskip.3cm

\hskip.4cm  i) $\beta$ induces, by restriction, an \'etale covering of $\bold P^1_K-\{0,1,\infty\}$.

\hskip.4cm ii) $\beta$ is tamely ramified above $\{0,1,\infty\}$ and totally ramified above $\infty$.

\hskip.4cm iii) The fiber $\beta^{-1}(0)$ consists of $n$ points having $a_1,\dots,a_n$ as ramification indices. 

\hskip.4cm iv) There is only one ramified point $x_\beta$ in the fiber $\beta^{-1}(1)$.

\vskip.3cm
  
\noindent In this case, $\te{deg}(\beta)=a_1+\dots+a_n$ and the Riemann-Hurwitz formula implies that $n$ is the ramification index of $\beta$ at $x_\beta$. Two such coverings $\beta_1$ and $\beta_2$ are isomorphic if and only if there exists an automorphism $\phi\in\te{PGL}_2(K)$ of the projective line over $K$ such that $\beta_2=\beta_1\circ\phi$. By assumption, the ramification is tame, so that  the characteristic of $K$ does not divide the integer $na_1\dots a_n(a_1+\dots+a_n)$. There exist finitely many diameter four trees of given type (over $K$). 
In characteristic zero, they can easily be enumerated using  the Grothendieck's theory of dessins d'enfants. More precisely,  when applied to coverings, the topological theory of the fundamental group induces a bijection between the set of diameter four trees of type $(a_1,\dots,a_n)$  and the set of (usual combinatorial) planar trees  endowed with a bipartite structure (a distinction between white and black vertices in such a way that the two ends of any edge never have the same color) having $n$ black vertices, of valencies $a_1,\dots,a_n$  and only one white vertex of valency greater than one (and then equal to $n$). The following picture is an illustration of this construction for diameter four trees of type (1,2,3) and (1,1,2,3). 

\vskip.8cm

\epsfysize=7.8cm

\centerline{\epsfbox{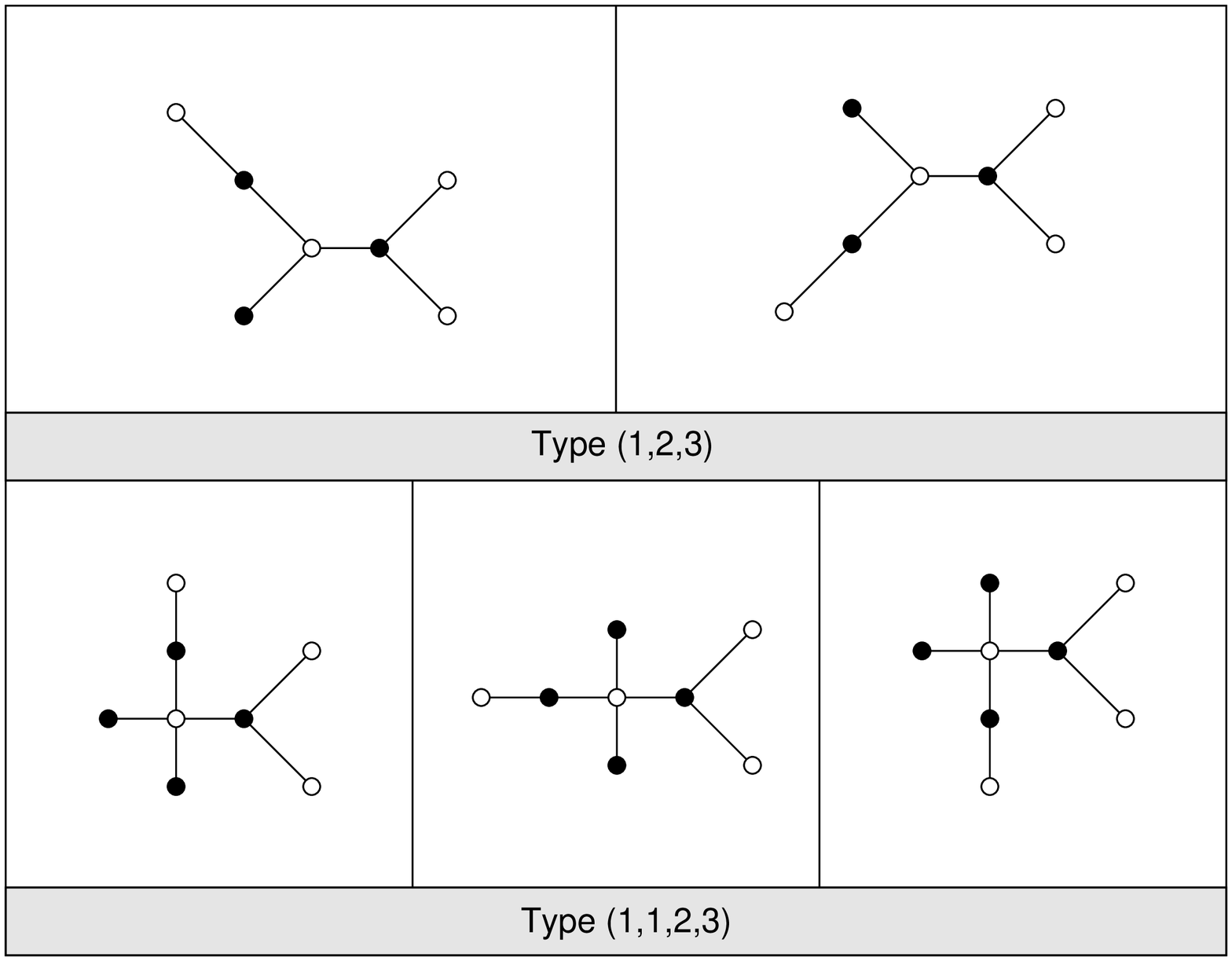}}

\vskip.6cm

 In particular, in the {\bbol generic} case, i.e., if $a_1<\dots<a_n$, there exist exactly $(n-1)!$ diameter four trees of type $(a_1,\dots,a_n)$ over $\QQ$. We refer to [G] for the first historical introduction to dessins d'enfants. A more detailed general treatement can be found in  [SV] and in the collected papers [S]. For further arithmetic topics in the genus zero case, see [K], [Z1] and [Z2]. 

\vskip1.2cm

\specialhead\twoplverm 2 Standard and normalized models\endspecialhead 

\vskip.3cm 

The notations being as in the previous paragraph, let $\Cal T=[\beta]$ be a diameter four tree over $K$, of type $(a_1\dots,a_n)$. The covering $\beta$ can be assimilated to an inclusion of fields $K(T)@>\quad>>K(X)$, sending $T$ to a rational function $\beta(X)\in K(X)$. Via the PGL$_2(K)$-action, we can always reduce to the case $\beta(\infty)=\infty$, so that $\beta(X)\in K[X]$ is a polynomial, since it only has one pole. Such a polynomial is not unique and will be called a (polynomial) {\bbol model} for the diameter four tree $\Cal T$. A model $\beta(X)$ is {\bbol standard} (resp. {\bbol normalized}) if $x_\beta=0$ (resp. $x_\beta=0$ and $\beta(1)=0$). For any model $\beta(X)$, the model $\beta(X+x_\beta)$ is standard. Two standard models $\beta_1(X)$ and $\beta_2(X)$ are associated to the same diameter four tree if and only if there exists and element $a\in K^*$ such that $\beta_2(X)=\beta_1(aX)$. A diameter four tree $\Cal T$ has finitely many normalized models. They can be constructued  as follows: let $\beta(X)$ be any standard model for $\Cal T$. Then, for any $x\in\beta^{-1}(0)$, the polynomial $\beta(x^{-1}X)$ is a normalized model associated to $\Cal T$, and they are all obtained in this way. A diameter four tree $\Cal T$ of type $(a_1,\dots,a_n)$   possesses $n\over m$ normalized models, where $m$ is the order of the automorphism group of $\Cal T$. 

\vskip1.2cm

\specialhead\twoplverm  3 Galois action, fields of moduli\endspecialhead 

\vskip.3cm

The normalized models associated to a diameter four tree $\Cal T$ over $K$ are automatically defined over a finite extension of the prime field $F$ of $K$. Moreover, if $\beta_1(X),\beta_2(X)\in\overline F[X]$ are two models associated to $\Cal T$, i.e., if $\beta_2(X)=\beta_1(aX+b)$ (with $a\in K^*$ and $b\in K$), then we have $a\in\overline F^*$ and $b\in\overline F$. For these reasons, we can reduce to the case  $K=\QQ$ or $K=\Fp$, where $p$ is a prime number. We then get a natural action of the absolute Galois group $G_F=\te{Gal}(K/F)$ on the set of  (normalized) models, which is compatible with the notion of isomorphism and induces a well defined  action of $G_F$ on the set of diameter four trees over $K$.  If $\beta(X)\in K[X]$ is a model associated to a tree $\Cal T$, then the {\bbol field of moduli} of $\Cal T$ is the subfield $F(\Cal T)$ of $K$ fixed by the elements $\sigma\in G_F$ such that $^\sigma\beta\cong\beta$, i.e., there exist elements $a_\sigma\in K^*$ and $b_\sigma\in K$ such that $^\sigma\beta(X)=\beta(a_\sigma X+b_\sigma)$. This field only depends on $\Cal T$ and not on the particular model. It is the intersection of the fields of definition of all the models associated to $\Cal T$. The index $[F(\Cal T):F]$ is the cardinality of the galois orbit of $\Cal T$. Moreover, using Hilbert '90, we can easily construct a standard model for $\Cal T$ defined  over $F(\Cal T)$. If the diameter four tree  is generic, then its field of moduli is the field of definition of any normalized model. In general, if $\beta(X)$ is a normalized model associated to $\Cal T$, the {\bbol splitting field} $F(\Cal T)^\bullet$ is the extension of $F$ obtained by adjunction of the elements belonging to the fiber $\beta^{-1}(0)$. It is easily seen that this field only depends on the tree and not on the particular model. It is a finite Galois extension of $F(\Cal T)$. 

\vskip.1cm 

The set IV$_{a_1,\dots,a_n}(K)$ consisting of all the diameter four trees of type $(a_1,\dots,a_n)$ over $K$ is called a {\bbol valency class}. Since the ramification indices are Galois invariants, we see that valencies classes are stable under the action of $G_F$. The compositum $F(\te{IV}_{a_1,\dots,a_n})$ of the fields of moduli of all the diameter four trees of type $(a_1,\dots,a_n)$ over $K$ is a  finite Galois extension of $F$. The splitting field $F(\te{IV}_{a_1,\dots,a_n})^\bullet$ is defined in a similar way.  One of the central questions in the Grothendieck theory of dessins d'enfants is to determine how a given valency class decomposes in Galois orbits.

\vskip1.2cm

\specialhead\twoplverm 4 Fundamental equations for standard models\endspecialhead 

\vskip.3cm 

\es The goal of this section is to give an algebraic characterisation of standard models associated to  diameter four trees. In the following, if $k\leq m$ are two non negative integers, then $\binom mk={m!\over k!(m-k)!}$ will denote the usual binomial coefficient. For simplicity, we have set $\binom mk=0$ for $k>m$. 

\vskip.4cm 

\proclaim{\twolverm 4.1 Proposition}\ita Let $K$ be an algebraically closed field and consider $n$ positive integers $a_1\leq\dots\leq a_n$. Suppose that the  characteristic of $K$ does not divide  $na_1\dots a_n(a_1+\dots+a_n)$. Let $x_1,\dots,x_n$ be pairwise ditinct elements of $K^*$ and set $\beta(X)=\prod_{i=1}^n(1-x_iX)^{a_i}\in K[X]$. Consider the following conditions:

\vskip.4cm

\noindent i) The polynomial $\beta(X)$ is a standard model associated to a diameter four tree of 

 type $(a_1,\dots,a_n)$. 

\vskip.2cm

\noindent ii) The elements $x_1,\dots,x_n$ define a solution of the system of algebraic equations  

$\psi_1=\dots=\psi_{n-1}=0$, where, for any $m\in\{1,\dots,n-1\}$,  

$$\psi_m=\sum_{k_1+\dots+k_n=m}\prod_{i=1}^n\binom{a_i}{k_i}x_i^{k_i}$$

\vskip.2cm 

\noindent iii) The elements $x_1,\dots,x_n$ define a solution of the system of algebraic equations  

$\phi_1=\dots=\phi_{n-1}=0$, where, for any $m\in\{1,\dots,n-1\}$,  

$$\phi_m=\sum_{i=1}^na_ix_i^m$$

\vskip.2cm

\noindent iv) For each $i\in\{1,\dots,n\}$, we have the identity 

$$a_i\prod_{j\neq i}(x_j-x_i)=(a_1+\dots+a_n)\prod_{j\neq i}x_j$$
 
\vskip.4cm 

\noindent With these notations and hypothesis, we have i) $\,\,\Leftrightarrow\,\,$ ii) $\,\,\Rightarrow\,\,$ iii) $\,\,\Leftrightarrow\,\,$ iv). Furthermore, if the  characteristic of $K$ does not divide $(n-1)!$, then all these conditions are equivalent.

\endproclaim

\vskip.4cm

\demo{\twolverm Proof}\es  We will start by proving the equivalence of conditions i) and ii). Denote by $\nu_0:K(X)^*@>\quad>>\bold Z$ the valuation associated to the point $X=0$. Then, the polynomial $\beta(X)=\prod_{i=1}^n(1-x_iX)^{a_i}\in K[X]$
is a standard model associated to a diameter four tree of type $(a_1,\dots,a_n)$ if and only if $\nu_0(\beta-1)=n$. If $\beta(X)=1+c_1X+\dots+c_NX^N$ (with $N=a_1+\dots+a_n$), then this last condition can be restated as $c_1=\dots=c_{n-1}=0$. Now, the expression of $\beta(X)$ in terms of $x_1,\dots,x_n$ gives 
$$\beta(X)=\prod_{i=1}^n(1-x_iX)^{a_i}=\prod_{i=1}^n\left(\sum_{k=1}^{a_i}(-1)^k\binom {a_i}kx_i^kX^k\right)=$$

$$=\sum_{m=0}^n(-1)^mX^m\sum_{k_1+\dots+k_n=m}\prod_{i=1}^n\binom{a_i}{k_i}x_i^{k_i}$$
In particular, we obtain $c_m=(-1)^m\psi_m$, and deduce from this the equivalence of the first two conditions. 

Since the characteristic $p$ of $K$ does not divide $n$, the condition $\nu_0(\beta-1)=n$ implies $\nu_0(\beta'/\beta)=n-1$, the converse being true if $p$ does not divide $(n-1)!$. This last condition can be handled more easily, and we will show that it is equivalent to the relations in iii) and iv). Let $\Cal O$ be the local ring at the point $X=0$, i.e., $\Cal O$ is the ring of integers of the valuation $\nu_0$. Its maximal ideal $\frak m$ is generated by the element $X\in K[X]$, its $\frak m$-adic completion $\widehat{\Cal O}$ is isomorphic to the formal power series ring $K[\![X]\!]$ and the canonical map $\Cal O@>\quad>>\widehat{\Cal O}$ is injective.  In the ring $\widehat{\Cal O}$, we have the identities
$${\beta'(X)\over\beta(X)}=-\sum_{i=1}^n{a_ix_i\over 1-x_iX}\quad\te{and}\quad{1\over1-x_iX}=\sum_{k\geq 0}x_i^kX^k,$$
so that the coefficient of $X^k$ in the power series expansion of $\beta'/\beta$ is $-\phi_{k+1}$.  The condition $\nu_0(\beta'/\beta)=n-1$ is then equivalent to the vanishing of these coefficients for $k=0,\dots,n-2$, which leads to the system in condition iii). In order to conclude the proof, we just have to show that conditions iii) and iv) are equivalent. Assume that iii) holds. Remark that the system of equations  is linear on $a_1,\dots,a_n$. If we consider the further relation  $N=a_1+\dots+a_n$, we then get a system of $n$ linear equations in $n$ variables, and a direct computation leads to the relations in iv), which are the expressions of $a_1,\dots,a_n$, in terms of $x_1,\dots,x_n$ and $N$. Finally, suppose that condition iv) holds.   Let's start with a general consideration: let $x_1,\dots,x_n$ be pairwise distinct elements of $K^*$.  For any $i\in\{1,\dots,n\}$, set $y_i=\prod_{j\neq i}(x_i-x_j)^{-1}$. We then have the following identity, which is easily proved:
$$\sum_{i=1}^n{y_i\over 1-x_iX}={X^{n-1}\over\prod_{i=1}^n(1-x_iX)}$$
In the present situation, the equations in iv) can be rewritten as $$a_ix_i=(-1)^{n-1}Nx_1\dots x_ny_i$$
so that, if we  multiply the above relation  by $-x_1\dots x_n(a_1+\dots+a_n)$, we  get 
$$-\sum_{i=1}^n{a_ix_i\over 1-x_iX}=(-1)^nx_1\dots x_n(a_1+\dots+a_n){X^{n-1}\over\prod_{i=1}^n(1-x_iX)}$$
Since the left side of this equality is $\beta'/\beta$, whe obtain $\nu_0(\beta'/\beta)=n-1$, which is equivalent to  condition ii). This concludes the proof of the proposition.  Remark that these last identities lead to the  expression
$$\beta'(X)=(-1)^nx_1\dots x_n(a_1+\dots+a_n)X^{n-1}\prod_{i=1}^n(1-x_iX)^{a_i-1}$$
which will be very usefull later on.\qed\enddemo

\vskip.4cm

\example{\twolverm 4.3 Examples} \es  

\vskip.2cm 

{\bbol 4.3.1 Diameter four trees of type $\bold {(a,b)}$ and $\bold {(a,b,c)}$.} We will now give a complete description of diameter four trees of type $(a,b)$, and $(a,b,c)$, for which the equations of proposition 4.1 can be easily solved.   The following picture is a combinatorial description of these diameter four trees over $\QQ$.

\vskip.6cm

\epsfysize=4cm

\centerline{\epsfbox{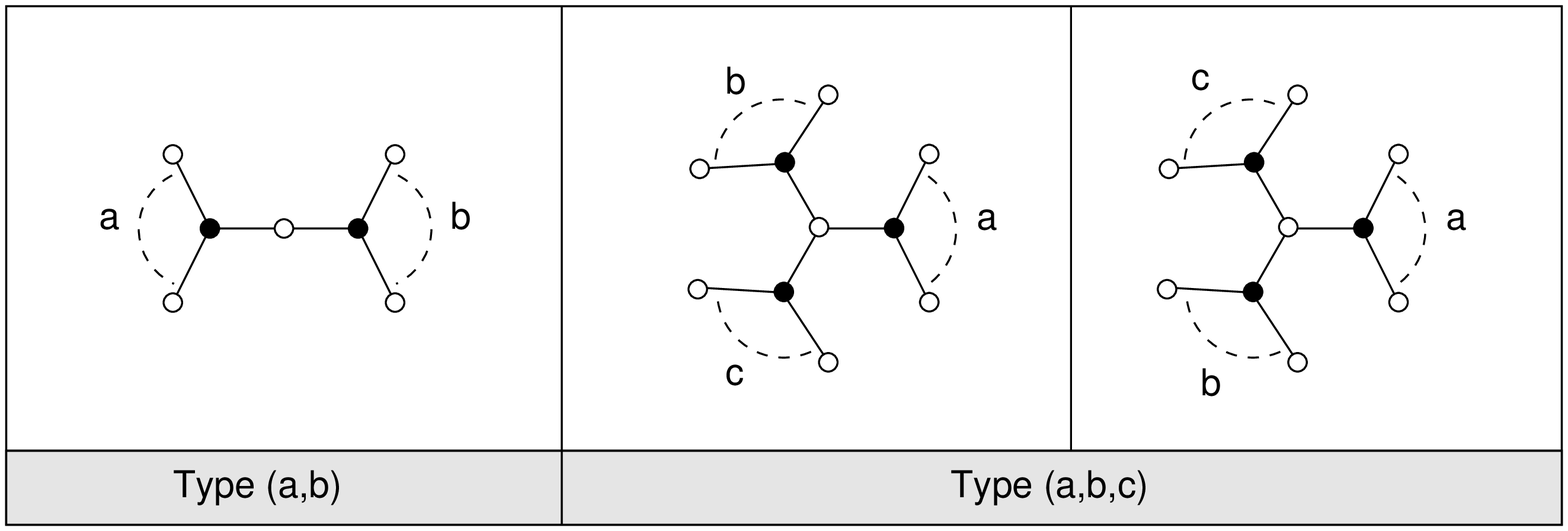}}

\vskip.6cm 

We see that the valency class IV$_{a,b}(\QQ)$, contains a unique element $\Cal T$, which is automatically defined over $\bold Q$. The same situation holds for non generic diameter four trees of type $(a,b,c)$. In the generic case,  there exist exactly two distinct diameter four trees in IV$_{a,b,c}(\QQ)$, but this description does not allow us to predict if they are conjugated or not under the action of $G_\bold Q$. In order to get some more subtle arithmetic informations, we will use the equations of the proposition. For trees of type $(a,b)$, we obtain the single equation $ax+by=0$, and  the polynomial $\beta(X)=(1-bX)^a(1-aX)^b$ is a standard model in any characteristic (not dividing $2ab(a+b)$) having $F$ as splitting field. We then deduce that IV$_{a,b}(K)$ conatains a unique element for any field $K$. In the case of diameter four trees of type $(a,b,c)$, and for $p>3$,  we have to solve  the system of algebraic equations
$$\left\{\aligned 
& ax+by+cz=0 \\
& ax^2+by^2+cz^2=0 \endaligned\right.$$
where we can set $z=1$ (and then obtain a normalized model). It turns out that over $\QQ$, in the non generic case, the only diameter four tree $\Cal T$ of type $(a,b,c)$ satisfies $\bold Q(\Cal T)=\bold Q$ and $\bold Q(\Cal T)^\bullet=\bold Q(\sqrt{t})\neq\QQ$, where $t^2+abc(a+b+c)=0$. In the generic case, we find that the field of moduli (and also the splitting field) of the two diameter four trees of type $(a,b,c)$ is the imaginary quadratic field $\bold Q(\sqrt{t})$, where $t$ is defined as above. In particular, they are conjugated, so that IV$_{a,b,c}(\QQ)$ is always a Galois orbit.

Let's now describe the valency class IV$_{a,b,c}(\Fp)$, for $p$ not dividing $6abc(a+b+c)$. The  situation is more complicated and we have to distinguish between different cases. First of all, set $D=(a+b)(b+c)(c+a)$ and $d=(a+b,b+c)(b+c,c+a)(c+a,a+b)$, where $(r,s)$ denotes the greatest common divisor of $r$ and $s$. Here are the different possibilities: 

i) If $p$ does not divide $D$, then the system can be solved as in characteristic zero: in the non generic case, we find one diameter four tree, defined over  $\bold F_p$, and its splitting field is $\bold F_p$ if and only if $-abc(a+b+c)$ is a square (take for example $a=2$, $b=3$, $c=4$ and $p=11$). In the generic case, there are again two diameter four trees of type $(a,b,c)$ over $\Fp$ and their field of moduli is $\bold F_p(\sqrt{t})$, where $t$ is defined as above. In particular, IV$_{a,b,c}(\Fp)$ is a Galois orbit if and only if $-abc(a+b+c)$ is not a square in $\bold F_p$. 
 
ii) If $p$ divides $d$, then IV$_{a,b,c}(\Fp)=\emptyset$. Remark that in the non generic case, if $p$ divides $D$, then it automatically divides $d$.

iii) Suppose now that $p$ divides $D$ but not $d$. In this case, the diameter four trees under consideration are  automatically generic. After some easy calculations, we deduce  that IV$_{a,b,c}(\Fp)$ contains a unique element $\Cal T$, which satisfies $\bold F_p(\Cal T)=\bold F_p(\Cal T)^\bullet=\bold F_p$.

In order to have a complete description of diameter four trees of type $(a,b,c)$, we just need to study the case $p=2$. In this situation, the integers $a,b,c$ are odd, and we have to consider the system of equations of condition ii), which gives 
$$\left\{\aligned 
&x+y+z=0\\
&xy+yz+zx+{a(a-1)\over 2}x^2+{b(b-1)\over 2}y^2+{c(c-1)\over 2}z^2=0\endaligned\right.$$
where, as before, we can set $z=1$. We find that IV$_{a,b,c}(\overline{\bold F_2})\neq \emptyset$ if and only if $a,b$ and $c$ have the same residue modulo $4$ (i.e.,  they all are congruent to $1$ or $3$ modulo $4$). If this condition is fullfilled, then,  as in characteristic zero, in the non generic case we find only one diameter four tree, and in the generic case, there are two diameter four trees of type $(a,b,c)$, which are conjugated by $G_{\bold F_2}$. More precisely, if $x\in\bold F_4$ satisfies $x^2+x+1$, then   the polynomial $\beta(X)=(1-X)^a(1-xX)^b(1-(1+x)X)^c$ is a normalized model, and they are all obtained in this way.

\vskip.2cm

{\bbol 4.3.2 Diameter four trees of type $\bold{(1,\dots,1,a,b)}$.} We will conclude this section with an example in characteristic zero which will appear quite often in the rest of the paper. Consider two integers $a$ and $b$ such that $1<a<b$. We want to study diameter four trees of type $(a_1,\dotsm,a_n)$ over $\QQ$, where $a_1=\dots=a_{n-2}=1$, $a_{n-1}=a$ and $a_{n}=b$. The combinatorial description shows that there exist exactly $n-1$ diameter four trees on this valency class. The following picture illustrates the case $n=5$.

\vskip.6cm

\epsfysize=3.2cm

\centerline{\epsfbox{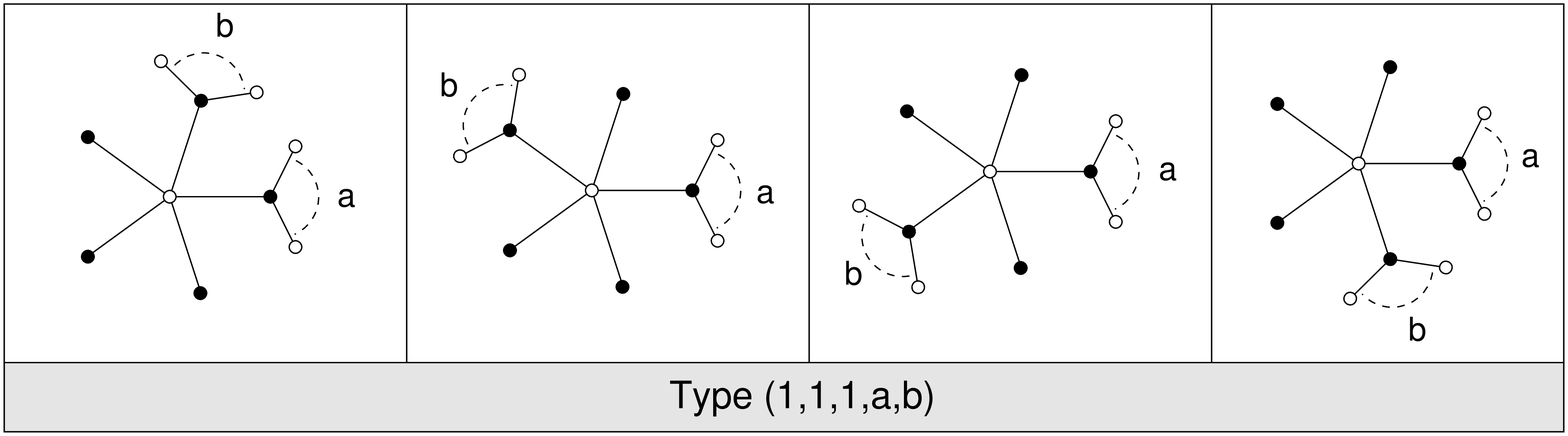}}

\vskip.6cm  If $\Cal T\in\te{IV}_{a_1,\dots,a_n}(\QQ)$, consider the (unique) normalized model $\beta(X)=\prod_{i=1}^n(1-x_iX)^{a_i}$ associated to $\Cal T$ such that  $x_n=1$ and set $x=x_{n-1}$. Then, it is easily shown that the field of moduli of $\Cal T$ is just $\bold Q(x)$. In order to solve the problem, we won't use the definig equations of proposition 4.1, but just remark that the derivative of $\beta(X)$ is given by $\beta'(X)=uX^{n-1}(1-xX)^{a-1}(1-X)^{b-1}$, with $u\in\QQ$. From this expression, we can deduce a formula for $\beta(X)$, only depending on $a,b,n$ and $x$. Moreover, the relation $\beta(1)=\beta(x)$ allows us to find all the possible values of $x$, and thus to completely solve the problem. After some elementary semplifications, following from the properties of binomial coefficients, we find that $x$ is a root of the polynomial 
$$h(X)=\sum_{k=0}^{n-1}\binom{a+k-1}{a-1}\binom{b+n-2-k}{b-1}X^k$$
For example, if we set $n=3$, we then obtain $h(X)={b(b+1)\over 2}+abX+{a(a+1)\over 2}X^2$, so that $\bold Q(x)=\bold Q(\sqrt{-ab(a+b+1)})$, which agrees with the results of the previous example. We have to remark that, even if we have an explicit expression for the field of moduli, it is not easy to determine whether or not IV$_{1,\dots,1,a,b}(\QQ)$ is a Galois orbit.  We already know that this is the case for $n=3$, and many direct calculations suggest that  this valency class is actually a Galois orbit for any $a,b$ and $n$. As we will see in the next sections, we will be able to prove this assertion in many cases, but in the more general setting it still remains a conjecture.

\vskip.2cm 

{\bbol 4.3.3 Diameter four trees with $\bold F_p$ as splitting field.}   We will now treat an example in positive characteristic $p$, assuming  that $p>n$. Let $u\in\bold F_p^*$ and consider $n$ pairwise distinct elements $x_1,\dots,x_n\in\bold F_p^*$ (they exist, since we are supposing $p>n$). For each $i\in\{1,\dots,n\}$, fix a positive integer $a_i$ whose reduction modulo $p$ is $u\prod_{j\neq i}(1-x_ix_j^{-1})^{-1}$. Then, the condition iv) of the proposition implies that the polynomial $\beta(X)=\prod_{i=1}^n(1-x_iX)^{a_i}$ is a standard model associated to a diameter four tree $\Cal T$ of type $(a_1,\dots,a_n)$ (up to reordering) such that  $\te{deg}(\beta)\equiv u\,\,\te{mod}(p)$  and $\bold F_p(\Cal T)^\bullet=\bold F_p$. Remark that we may choose $a_1,\dots,a_n$ in such  a way that the diameter four tree $\Cal T$ is generic. As a numerical example, set $p=11$, $n=5$, $x_1=1, x_2=2, x_3=3, x_4=5, x_5=6$ and $u=9$. Then, the polynomial 
$$\beta(X)=(1-3X)(1-X)^2(1-5X)^3(1-2X)^4(1-6X)^{10}$$
is a normalized model for a diameter four tree of type $(1,2,3,4,10)$ over $\overline{\bold F}_{11}$.

\endexample

\vskip1.2cm

\specialhead{\twoplverm 5 p-congruence}\endspecialhead

\vskip.3cm

This section is devoted to the study of a general construction, which only works in positive characteristic. Its applications will be very important when applied to the lifting of normalized models in characteristic zero and, in many cases,  will allow us to completely  describe the action of decomposition groups. The starting point is the following elementary lemma:

\vskip.4cm

\proclaim{\twolverm 5.1 Lemma}\ita Let $K$ be an algebraically closed field of characteristic $p>0$. Denote by  $\nu_0:K(X)^*@>\quad>>\bold Z$  the valuation associated to the point $X=0$. Let $\beta(X)\in K[X]$ be a polynomial such that $\beta(0)=1$. Let $h,n$ be positive integers such that $p^h>n$ and consider an element $x\in K^*$. Then, for any positive integer $m$, setting $\beta_1(X)=(1-xX)^{mp^h}\beta(X)$, we have $\nu_0(\beta-1)\geq n$ if and only if $\nu_0(\beta_1-1)\geq n$.  
\endproclaim

\vskip.4cm

\demo{\twolverm Proof}\es Set $\beta(X)=1+X^kq(X)$, with $q(X)\in K[X]$ and $q(0)\neq 0$. Then,  $\beta_1(X)=(1-xX)^{mp^h}\beta(X)=(1-(xX)^{p^h})^m\beta(X)=(1-X^{p^h}t(X))(X^kq(X)+1)=1+X^kq(X)+X^{p^h}s(X)$, with $s(X),t(X)\in K[X]$. In particular, if $k<p^h$ then  $\nu_0(\beta-1)=k=\nu_0(\beta_1-1)$, and the lemma follows.\qed
\enddemo

\es We want to apply this result to diameter four trees. In general, if $a_1\leq\dots\leq a_n$ and $b_1\leq\dots\leq b_n$  are positive integers, we will say that $(a_1,\dots,a_n)$ and $(b_1,\dots,b_n)$ are {\bbol $\bold p$-congruent} (resp. {\bbol strictly $\bold p$-congruent}) if there exists a permutation $\sigma\in\frak S_n$ such that $a_i\equiv b_{\sigma(i)}\,\,\te{mod}(p^h)$ for any $i\in\{1,\dots,n\}$ (resp. $a_i\equiv b_{\sigma(i)}\,\,\te{mod}(p^h)$ and $a_i=a_j$ if and only if $b_{\sigma(i)}=b_{\sigma(j)}$), where $h=h_p(n)$ is the least integer such that $p^h>n$. The permutation $\sigma$ is clearly not unique; we will say that it is {\bbol admissible} (with respect to $(a_1,\dots,a_n)$ and $(b_1,\dots,b_n)$). 

\vskip.6cm 

\proclaim{\twolverm 5.2 Proposition} Let $\beta(X)=\prod_{i=1}^n(1-x_iX)^{a_i}$ be a standard model associated to a diameter four tree $\Cal T$ of type $(a_1,\dots,a_n)$ over a field $K$ of positive characteristic $p>0$. Suppose that $(b_1,\dots,b_n)$ is $p$-congruent to $(a_1,\dots,a_n)$ and let $\sigma$ be an admissible permutation. Then, the polynomial $\beta_1(X)=\prod_{i=1}^n(1-x_iX)^{b_{\sigma(i)}}$ is a standard model associated to a diameter four tree $\Cal T_1$ of type $(b_1,\dots,b_n)$. \endproclaim

\vskip.4cm

\demo{\twolverm Proof}\es Without any loss of generality, we can reduce to the case $b_{\sigma(1)}=a_1+mp^h$, with $m>0$, and $b_{\sigma(i)}=a_i$ for $i>1$. We then have $\beta_1(X)=(1-x_1X)^{p^hm}\beta(X)$, and the corollary follows from lemma 5.1 and from the fact that $\beta_1(X)$ is a standard model if and only if $\nu_0(\beta_1-1)=n$\qed.
\enddemo

\vskip.2cm

\es  Following the notations and hypothesis of this proposition, we have  $\beta^{-1}_1(0)=\beta^{-1}(0)$, so that $\bold F_p(\Cal T_1)^\bullet=\bold F_p(\Cal T)^\bullet$. As a direct consequence, we get $\bold F_p(\te{IV}_{a_1,\dots,a_n})^\bullet=\bold F_p(\te{IV}_{b_1,\dots,b_n})^\bullet$. The definition of $\beta_1(X)$ clearly depends on the admissible permutation $\sigma$. This problem disappears if we suppose that  $(a_1,\dots,a_n)$ is strictly $p$-congruent to $(b_1,\dots,b_n)$. Indeed, we have the following straightforward result:

\vskip.4cm

\proclaim{\twolverm 5.3 Proposition}\ita Let $a_1\leq\dots\leq a_n$ and $b_1\leq\dots\leq b_n$ be positive integers and consider a prime number $p$. If $(a_1,\dots,a_n)$ is strictly $p$-congruent to $(b_1,\dots,b_n)$ then there is a bijection between IV$_{a_1,\dots,a_n}(\Fp)$ and IV$_{b_1,\dots,b_n}(\Fp)$ which commutes with the action of $G_{\bold F_p}$. In particular, we have  $\bold F_p(\te{IV}_{a_1,\dots,a_n})=\bold F_p(\te{IV}_{b_1,\dots,b_n})$.\endproclaim

\vskip.4cm

\es  Finally, these constructions show that when studying (the Galois action on) diameter four trees of type $(a_1,\dots,a_n)$ in positive characteristic $p$, one can always reduce to the case $0<a_1\leq\dots\leq a_n\leq n(p-1)$ (in the generic case, we can even suppose $0<a_1\leq\dots\leq a_n<p$).

\vskip.4cm

\example{\twolverm 5.4 Example}\es Let $a_1\leq\dots\leq a_n$ be positive integers such that $(a_1,\dots,a_n)$ is $p$-congruent to $(1,\dots,1)$, i.e.,  for any $i\in\{1,\dots,n\}$, we have $a_i\equiv 1\,\,\te{mod}(p^h)$, where $h=h_p(n)$. Suppose that $\beta(X)=\prod_{i=1}^n(1-x_iX)^{a_i}$ is a normalized model for a diameter four tree of type $(a_1,\dots,a_n)$ over $\Fp$. Then, the polynomial $\delta(X)=\prod_{i=1}^n(1-x_iX)$ is a normalized model for a diameter four tree of type $(1,....,1)$. Now, since $p$ does not divide $n$, we see that IV$_{1,\dots,1}(\Fp)$ consists of only one element, having the polynomial $\rho(X)=1-x^n$ as unique normalized model. We then deduce that $\delta(X)=\rho(X)$. In particular, the elements $x_1,\dots,x_k$ all belong to the group of $\mu_n$ of $n$-th roots of unity. The converse is true. Namely, given pairwise different (ordered) elements $x_1,\dots,x_n\in\mu_n$, the polynomial $\beta(X)$ defined as above is a normalized model associated to a diameter four tree of type $(a_1,\dots,a_n)$. We deduce from this that in the generic case, the cardinality of IV$_{a_1,\dots,a_n}(\Fp)$ is $(n-1)!$ and decomposes in $(n-1)!\over\te{\esp{ord}}_n(p)$ Galois orbits, where $\te{ord}_n(p)$ is the order of $p$ in $(\bold Z/n\bold Z)^*$ (a similar, but more complicated formula can be obtained in the non generic case). 

\endexample

\vskip1.2cm

\specialhead{\twoplverm 6  Lifting normalized models from positive characteristic}\endspecialhead

\vskip.3cm

We could now adopt two different points of view in order to study the problem of lifting normalized models from positive characteristic: first of all, we could interpretate the equations of proposition 4.1 in a purely algebraic-geometrical way, by constructing a (projective) scheme whose points naturally correspond to diameter four trees, and then deduce some arithmetical results arising from the general machinery of the theory. Another possible approach is to work directly with models and then only use a few fundamental results, as the (multidimensional) Hensel's lemma. Clearly, the first point of view leads to finer results and is somewhat more natural, but it may hide  the philosophy and the essential steps of the construction. That's why we decided to study the problem in the more elementary way, and thus adopt the second point of view.

\vskip.4cm

\proclaim{\twolverm 6.1 Proposition}\ita  Let $(R,\frak p)$ denote a complete discrete valuation ring of field  of fraction $K$ of characteristic $0$ and of residue field $k=R/\frak p$ of characteritic $p>0$. Suppose that $\overline\beta(X)\in k[X]$ is a normalized model associated to a diameter four tree of type $(a_1,\dots,a_n)$. Then, it can be uniquely lifted to a  normalized model  $\beta(X)\in R[X]$ for a diameter four tree  of the same type. \endproclaim 

\vskip.3cm

\demo{\twolverm Proof}\es Consider a Galois extension $R'$ of $R$ of field of fraction $K'$ and residue field $k'$, such that $\bold F_p(\Cal T_p)^\bullet\subset k'$.  
 Set $\overline\beta(X)=\prod_{i=1}^n(1-\overline x_iX)^{a_i}$, with $\overline x_1,\dots,\overline x_n\in (k')^*$ pairwise distinct and one of them equal to $1$. We may suppose, without loss of generality, that $\overline x_n=1$. 
Following proposition 4.1, the  elements $\overline x_1,\dots,\overline x_{n-1}$ define a solution of  the system of $n-1$ algebraic equations $\psi_1=\dots=\psi_{n-1}=0$, where $\psi_m=\sum_{k_1+\dots+k_n=m}\prod_{i=1}^n\binom{a_i}{k_i}X_i^{k_i}$ and $X_n=1$. The multidimensional Hensel's lemma then implies that if the determinant of jacobian matrix $J_\psi=(\partial\psi_i/\partial X_j)$ is a unit in $k'$ for $X_1=\overline x_1,\dots, X_n=\overline x_n=1$,  then this solution  can be uniquely lifted to a solution $x_1,\dots,x_{n-1},x_n=1$ with $x_i\in R'$ for any $i\in\{1,\dots,n-1\}$. A tedious but straightforward calculation shows that det$(J_\psi)={1\over(n-1)!}\te{det}(J_\phi)$, where $J_\phi=(\partial\phi_i/\partial X_j)$ is the jacobian matrix associated to the system of equations $\phi_1=\dots=\phi_n=0$, with $\phi_m=\sum_{i=1}^ma_ix_i^m$. Now, we have $\partial\phi_i/\partial X_j=ia_jX_j^{i-1}$, so that 
$$\te{det}(J_\phi)=\te{det}\left(\matrix a_1 & \dots  & a_{n-1}\\ \vdots & \, & \vdots\\ (n-1)a_1X_1^{n-2}&\dots & (n-1)a_{n-1}X_{n-1}^{n-2}\endmatrix\right)=$$
$$=(n-1)!a_1\dots a_{n-1}\te{det}\left(\matrix 1 & \dots  & 1\\ \vdots & \, & \vdots\\ X_1^{n-2}&\dots & X_{n-1}^{n-2}\endmatrix\right)=$$
$$=(n-1)!a_1\dots a_{n-1}\prod_{1\leq i<j\leq n-1}(X_i-X_j),$$
and we finally get 
$$\te{det}(J_\psi)=a_1\dots a_{n-1}\prod_{1\leq i<j\leq n-1}(X_i-X_j)$$
which is different from zero for $X_1=\overline x_1,\dots,X_n=\overline x_n=1$, since $\overline x_i\neq\overline x_j$ for any $i\neq j$ (by hypothesis) and  $p$ does not divide  $a_1\dots a_n$ (the ramification being tame). The corresponding normalized model in characteristic zero is given by $\beta(X)=\prod_{i=1}^n(1-x_iX)^{a_i}$. Let now $\sigma\in\te{Gal}(K'/K)=\te{Gal}(k'/k)$. The polynomial $^\sigma\beta(X)$ obtained by acting on the coefficients of $\beta(X)$ reduce to $\overline\beta(X)$ modulo $\frak p$ (since $\overline\beta(X)\in k[X]$). By uniqueness of the lifting, we then deduce that $^\sigma\beta(X)=\beta(X)$, so that $\beta(X)\in R[X]$.  
\qed  
\enddemo

\vskip1.2cm

\specialhead{\twoplverm 7 Good reduction, action of the decomposition groups}\endspecialhead

\vskip.3cm

The lifting result of the previous paragraph has many arithmetical applications.  We will start by introducing the notion of good reduction. First of all, let's fix some notations: from now on $K$  will denote a finite Galois extension of $\bold Q$ of group $G$. Let $\Cal O_\frak p$ be the localisation of its ring of integers $\Cal O_K$  at a (maximal) prime ideal $\frak p$. It is a discrete valuation ring. Let moreover $p$ be the prime number such that $\frak p\cap\bold Z=p\bold Z$ and denote by $\widehat{\Cal O}_\frak p$ the $\frak p$-adic completion of $\Cal O_\frak p$, and by $K_\frak p$ the field of fractions of $\widehat{\Cal O}_\frak p$.  The residue field $k(\frak p)=\Cal O_\frak p/\frak p=\widehat{\Cal O}_\frak p/\frak p\widehat{\Cal O}_\frak p$ is a finite Galois extension of the prime field $\bold F_p$, of group $D_\frak p/I_\frak p$ where $D_\frak p$ (resp. $I_\frak p$) is the decomposition group (resp. the inertia group) at $\frak p$. Moreover, we know that $\te{Gal}(k(\frak p)/\bold F_p)=\te{Gal}(K_\frak p/\bold Q_p)$. Finally, let $\pi$ be a uniformizer for $\frak p$ and denote by $\nu:K^*@>\quad>>\bold Z$ the associated valuation, normalized by $\nu(\pi)=1$.

\vskip.2cm

Let now $\beta(X)\in K[X]$ be a normalized model associated to a diameter four tree $\Cal T$ of type $(a_1,\dots,a_n)$.
We will say that $\Cal T$  has {\bbol good reduction} at $\frak p$ if there exists a normalized model $\beta(X)\in\Cal O_\frak p[X]$ such that its canonical image $\overline\beta(X)\in k(\frak p)[X]$ is a normalized model for a diameter four tree of the same type over $\Fp$. It is easily checked that if $\Cal T$ has good reduction at $\frak p$, then any (and not only one) normalized model satisfies this last condition. More precisely, we have the following

\vskip.4cm

\proclaim{\twolverm 7.1 Lemma}\ita With the above notations and hypothesis, suppose that $\bold Q(\Cal T)^\bullet\subset K$ (we may set $K=\bold Q(\te{IV}_{a_1,\dots,a_n})^\bullet$) and that $p$ does not divide $na_1\dots a_n(a_1+\dots+a_n)$. Then, the following conditions are equivalent:

\vskip.2cm

\noindent i) The diameter four tree $\Cal T$ has good reduction at $\frak p$.

\vskip.2cm 

\noindent ii) All the normalized models of $\Cal T$ are defined over $\Cal O_\frak p$.

\endproclaim

\vskip.4cm

\demo{\twolverm Proof}\es  i) $\Rightarrow$ ii) Let $\beta(X)=\prod_{i=1}^n(1-x_iX)^{a_i}\in\Cal O_\frak p[X]$ be a normalized model for $\Cal T$ such that  its reduction modulo $\frak p$ is a normalized model for a diameter four tree of the same type over $\Fp$. In particular, this implies  that $\overline x_1,\dots,\overline x_n\in k(\frak p)^*$, i.e., $x_1,\dots,x_n\in\Cal O_\frak p^*$. Now, all the other normalized models of $\Cal T$ are given by $\beta_i(X)=\beta(x_i^{-1}X)$, with $i\in\{1,\dots,n\}$, and since $x_i$ is a unit in $\Cal O_\frak p$, we deduce that $\beta_i(X)\in\Cal O_\frak p[X]$.

ii) $\Rightarrow$ i) Consider a normalized model  $\beta(X)=\prod_{i=1}^n(1-x_iX)^{a_i}\in\Cal O_\frak p[X]$ of $\Cal T$. Then, the polynomial $X^{a_1+\dots+a_n}\beta(X^{-1})=\prod_{i=1}^n(X-x_i)^{a_i}$ is monic and defined over $\Cal O_\frak p$. This implies that $x_1,\dots,x_n$ are integral over $\Cal O_\frak p$, which is integrally closed. Since we are assuming $\bold Q(\Cal T)^\bullet\subset K$, we obtain $x_1,\dots,x_n\in\Cal O_\frak p$. This holds for any normalized model, so that $x_ix_j^{-1}\in\Cal O_\frak p$ for any $i,j\in\{1,\dots,n\}$. Now, there is an $i\in\{1,\dots,n\}$ such that $x_i=1$, which implies that $x_1,\dots,x_n\in\Cal O_\frak p^*$. In order to prove that $\overline\beta(X)$ is a normalized model for a diameter four tree of the same type, we just have to check that $\overline x_i\neq\overline x_j$ for any $i\neq j$. This follows from the  identities in condition iv) of  proposition 4.1, since $p$ does not divide $a_1\dots a_n(a_1+\dots+a_n)$.\qed 

\enddemo

\vskip.4cm

\proclaim{\twolverm 7.2 Proposition}\ita Let $\Cal T$ be a diameter four tree of type $(a_1,\dots,a_n)$ such that $\bold Q(\Cal T)^\bullet\subset K$. Suppose that $\Cal T$ has good reduction at $\frak p$ and consider a normalized model $\beta(X)=\prod_{i=1}^n(1-x_iX)^{a_i}\in K[X]$ associated to  $\Cal T$. Then the inertia subgroup $I_\frak p$ acts trivially on $x_1,\dots,x_n$. In particular, if $\frak q=\frak p\cap\bold Q(\Cal T)^\bullet$, then the extension $\bold Q(\Cal T)^\bullet/\bold Q$ is unramified at $\frak q$. \endproclaim

\vskip.4cm

\demo{\twolverm Proof}\es  Denote by $\overline x_1,\dots,\overline x_n$ the canonical images of $x_1,\dots,x_n$ in $k(\frak p)$ and consider an element $\sigma\in I_\frak p$. By definition,  $\sigma$ acts trivially on $k(\frak p)$, so that $\overline{\sigma(x_i)}=\sigma(\overline x_i)=\overline x_i$ for any $i\in\{1,\dots,n\}$. In particular, $x_1,\dots,x_n$ and $\sigma(x_1),\dots,\sigma(x_n)$ define two solutions of the system in condition ii) of proposition 4.1 reducing to the same solution in characteristic $p$. Following proposition 6.1, $\overline x_1,\dots,\overline x_n$ can be uniquely lifted, and  we have $\sigma(x_i)=x_i$ for any $i\in\{1,\dots,n\}$, which completes the proof.\qed

\enddemo 

\vskip.4cm

We will now give a criterion for good reduction, only depending on the integers $a_1\leq\dots\leq a_n$. First of all, consider the integer $d(a_1,\dots,a_n)$ defined by  $$d(a_1,\dots,a_n)=\prod_{S\subset\{1,\dots,n\}}\sum_{i\in S}a_i$$
In other words, $d(a_1,\dots,a_n)$ is the product of the finite sums of elements of the set $\{a_1,\dots,a_n\}$.

\vskip.4cm

\proclaim{\twolverm 7.3 Proposition}\ita  Let $\Cal T$ be  a diameter four tree of type $(a_1,\dots,a_n)$ such that $\bold Q(\Cal T)^\bullet\subset K$. If $\frak p$ does not divide $d(a_1,\dots,a_n)$,  then $\Cal T$ has good reduction at $\frak p$. In particular, $\bold Q(\Cal T)^\bullet$ is unramified above the rational primes dividing $d(a_1,\dots,a_n)$.
\endproclaim

\vskip.4cm

\demo{\twolverm Proof}\es Let  $\beta(X)=\prod_{i=1}^n(1-x_iX)^{a_i}\in K[X]$ be a normalized model associated to $\Cal T$. We may suppose that $x_1,\dots,x_n\in\Cal O_\frak p$. Indeed, let $i\in\{1,\dots,n\}$ such that $\nu(x_i)=\te{Min}\{\nu(x_1),\dots,\nu(x_n)\}$. Then the polynomial $\beta(x_i^{-1}X)=\prod_{j=1}^n(1-y_jX)^{a_j}$ is a normalized model for $\Cal T$ satisfying the above condition. Suppose that $\frak p$ does not divide $d(a_1,\dots,a_n)$. In particular, $\frak p$ does not divide  $a_1\dots a_n(a_1+\dots+a_n)$.  Let $N_0=\{i\in\{1,\dots,n\}\,\,|\,\, \overline x_i\neq 0\}$. We have $N_0\neq\emptyset$, since the model is normalized (and thus, there exists $i\in\{1,\dots,n\}$ such that $x_i=1$). If $N_0=\{1,\dots,n\}$, then one can easily deduce from the relations in condition iv) of proposition 4.1  that $\overline x_i\neq\overline x_j$ for any $i\neq j$, so that $\Cal T$ has good reduction at $\frak p$.  Assume then that  $N_0$ is strictly contained in $\{1,\dots,n\}$. We can define an equivalence relation on $N_0$ by setting $i\sim j$ if and only if $\overline x_i=\overline x_j$. The quotient $N_1=N_0/\sim$ has cardinality $k<n$. For any $\tau\in N_1$, set $b_\tau=\sum_{i\in\tau} a_i$ and $x_\tau=\overline x_i$, with $i\in\tau$. Then, $x_\tau\neq x_\sigma$ for any $\tau\neq\sigma$ and the equations in condition iii) of proposition 4.1 give $\sum_{\tau\in N_1}b_\tau x_\tau^i=0$ for any $i\in\{1,\dots,n-1\}$. If we restrict to the first $k$ equations, we get a  homogeneous system of linear equations  on the $k$ variables $b_\tau$. Since its determinant does not vanish,  we then deduce that the only solution is given $b_\tau=0$ for any $\tau$, which can be rephrased as  $p$ divides $b_\tau$ for any $\tau$. In particular $p$ divides $d(a_1,\dots,a_n)$, which  contraddicts our assumption. We then have $N_0=\{1,\dots,n\}$, and  the diameter four tree $\Cal T$ has good reduction at $\frak p$, which concludes the proof.\qed 

\enddemo

\vskip.4cm 

We will now come ack to $p$-conruence and give some applications concerning the action of decomposition groups. First of all, for any diameter four tree $\Cal T$  over $\QQ$ such that $\bold Q(\Cal T)^\bullet\subset K$, set  $G_\Cal T=\te{Gal}(K/\bold Q(\Cal T))$ and  $G_\Cal T^\bullet=\te{Gal}(K/\bold Q(\Cal T)^\bullet)$. We have $G_\Cal T^\bullet\subset G_\Cal T\subset G$, and $G_\Cal T=G_\Cal T^\bullet$ in the generic case.

\vskip.4cm

\proclaim{\twolverm 7.4 Proposition}\ita Suppose that $(a_1,\dots,a_n)$ and $(b_1,\dots,b_n)$ are $p$-congruent. Assume that $K$ contains the compositum of $\bold Q(\te{IV}_{a_1,\dots,a_n})^\bullet$ and $\bold Q(\te{IV}_{b_1,\dots,b_n})^\bullet$. If $\Cal T$ is a  diameter four tree of type $(a_1,\dots,a_n)$ over $\QQ$ having good reduction at $\frak p$, then there exists a diameter four tree $\Cal T'$ of type $(b_1,\dots,b_n)$ over $\QQ$ having  good reduction at $\frak p$ such that  $G_\Cal T^\bullet\cap D_\frak p=G_{\Cal T'}^\bullet\cap D_\frak p$. \endproclaim 

\vskip.4cm

\demo{\twolverm Proof}\es It is a straightforward consequence of the uniqueness of the lifting (cf. proposition 6.1) and of the construction of \S 5.\qed
\enddemo

\vskip.4cm

In the case of strict $p$-congruence we can obtain some finer results, which directly follow from proposition 5.3:

\vskip.4cm

\proclaim{\twolverm 7.5 Proposition}\ita Under the same hypothesis (on $K$) of proposition 7.4, suppose now that $(a_1,\dots,a_n)$ and $(b_1,\dots,b_n)$ are strictly  $p$-congruent. If $\Cal T$ is a  diameter four tree of type $(a_1,\dots,a_n)$ over $\QQ$ having good reduction at $\frak p$, then there exists a canonical diameter four tree $\Cal T'$ of type $(b_1,\dots,b_n)$ over $\QQ$ having  good reduction at $\frak p$ such that  $G_\Cal T\cap D_\frak p=G_{\Cal T'}\cap D_\frak p$. \endproclaim 

\vskip.4cm

\example{\twolverm 7.6 Example}\es We will now obtain some arithmetical informations in characteristic zero, arising from the results of example 5.4 in positive characteristic.  Let $\Cal T$ be a diameter four tree of type $(a_1,\dots,a_n)$, with $(a_1,\dots,a_n)$ $p$-congruent to $(1,\dots,1)$ and suppose that $p$ does not divide $n$. After enlarging $K$, we may assume that $\bold Q(\Cal T)^\bullet\subset K$ and that $K$  contains the group  $\mu_n$ of $n$-th  roots of unity. Setting $G_n=\te{Gal}(K/\bold Q(\mu_n))$, we get an exact sequence 
$$1@>\quad>> G_n @>\quad>> G @>\,\,\,\,\chi\,\,\,\,>> (\bold Z/n\bold Z)^* @>\quad>>1$$
In example 5.3, we saw that in the generic case, the valency class IV$_{a_1,\dots,a_n}(\Fp)$ contains exactly $(n-1)!$ diameter four trees. In general, it can be easily shown that IV$_{a_1,\dots,a_n}(\Fp)$ and IV$_{a_1,\dots,a_n}(\QQ)$ have the same cardinalities  (always assuming that $(a_1,\dots,a_n)$ is  p-congruent to $(1,\dots,1)$). We deduce from this that $\Cal T$ has good reduction at any prime $\frak p$ lying over $p\bold Z$. There is only one diameter four tree $\Cal T_n$ of type $(1,\dots,1)$ over $\QQ$. One easily checks that $\bold Q(\Cal T_n)^\bullet=\bold Q(\mu_n)$, and the image of $D_\frak p$ in $(\bold Z/n\bold Z)^*$ is the subgroup generated by $p$. We then deduce from proposition 7.4 that the action of the decomposition group $D_\frak p$ on $\Cal T$ factors through $\bold Q(\mu_n)$, i.e., it only depends on the cyclotomic character $\chi$. 

\vskip.1cm 

As an application, let's come back to the diameter four trees introduced in example 4.3.2: suppose that $a_1=\dots=a_{n-2}=1$, $a_{n-1}=a$ and $a_n=b$, with $a<b$ and $a,b$ congruent to $1$ modulo $p^h$ (here, according with \S5, $h=h_p(n)$ is the least integer such that $p^h>n$). The valency class IV$_{a_1,\dots,a_n}(\QQ)$ contains $n-1$ diameter four trees. In terms of liftings from characteristic $p$, they can obtained in the following way: choose a $n$-th root of unity $\zeta\in\Fp$ and set $\overline\beta_\zeta(X)=(1-X)^{a-1}(1-\zeta X)^{b-1}\prod_{\zeta'\in\mu_n}(1-\zeta'X)$. Then, following proposition 6.1, the normalized model $\overline\beta_\zeta(X)$ can uniquely be lifted to a normalized model $\beta_\zeta(X)$. The  corresponding  diameter four tree will be denoted $\Cal T_\zeta$. Proposition 7.4 implies that the Galois orbit of $\Cal T_\zeta$ under the action of the decomposition group at $\frak p$ is the cardinality of the set $\{\zeta^{p^m}\,\,|\,\,m\in\bold Z\}$. Suppose, furthermore, that $n$ is a prime number and that $\te{ord}_n(p)=n-1$, i.e., that $p$ generates $\bold F_n^*$. In this case, we conclude that IV$_{1,\dots,1,a,b}(\QQ)$ is a Galois orbit and that, for any $\Cal T$ in this valency class, there is only one prime in $\bold Q(\Cal T)$ lying above $p\bold Z$. On the opposite situation, if $p$ is congruent to $1$ modulo $n$, then we cannot affirm that  IV$_{1,\dots,1,a,b}(\QQ)$ is a Galois orbit, but since in this case the decomposition groups act trivially, we deduce that $p\bold Z$ totally splits in the field of moduli of any such diameter four tree. 

As a numerical example, we can consider diameter four trees of type $(1,1,1,8a+1,8b+1)$, with $a<b$, and study the action of the decomposition groups of the primes lying over $2\bold Z$. In this case,  we have $h_2(5)=3$, so tat $(1,1,1,8a+1,8b+1)$ is $2$-congruent to $(1,1,1,1,1)$. Moreover, $n=5$ is prime and ord$_5(2)=4$. We deduce from this that IV$_{1,1,1,8a+1,8b+1}$ is always a  Galois orbit. Setting $a=1$ and $b=2$, we know from example 5.3, that $\bold Q(\te{IV}_{1,1,1,9,17})$ is the splitting field of the polynomial 
$$h(X)=495X^4+2805X^3+6885X^2+8721X+4845$$
It is a Galois extension of $\bold Q$ of group $\frak S_4$, and the reduction of $h(X)$ modulo $2$ is the irreducible cyclotomic polynomial $X^4+X^3+X^2+X+1$. Moreover, the odd primes dividing  $d(1,1,1,9,17)$ are $3,5,7,11,13,17,19$ and $29$, so that, following proposition 7.3, $\bold Q(\te{IV}_{1,1,1,9,17})$ is unramified above these primes. A direct computation of the discriminant gives $\Delta(h)=3^{10}\cdot5^2\cdot7^2\cdot11\cdot17^3\cdot19\cdot29^3$, and we see that the primes $11,17,19$ and $29$ effectively ramify (a finer analysis is needed for $3,5$ and $7$). 
\endexample

\vskip1.2cm

\specialhead{\twoplverm 8 Wild ramification above zero, Kummer models}\endspecialhead

\vskip.4cm

In the previous section, we studied the specialisation of diameter four trees  from characteristic zero to positive characteristic $p$, always assuming that $p$ does not divide $na_1\dots a_n(a_1+\dots+a_n)$, i.e., we restricted to primes not dividing the indices of ramification. We will now be concerned with the case of wild ramification above zero. More precisely, given $a_1\leq\dots\leq a_n$ and $i\in\{1,\dots,n\}$, set 
$$d_i(a_1,\dots,a_n)=d(a_1,\dots,\widehat a_i,\dots,a_n)=\prod_{S\subset\{1,\dots,\widehat i,\dots,n\}}\sum_{j\in S}a_j$$
Our purpose is to study the specialisation of diameter four trees at the primes dividing $a_i$ but not dividing $d_i(a_1,\dots,a_n)$ (for a fixed $i$). Such primes will be called {\bbol $\bold a_i$-regular}. Our first result is a description of the reduction of normalized models. As in the previous paragraph, $K$ will denote a number field, $\frak p$ a maximal ideal of $K$ and $p$ will be the prime number such that $\frak p\cap\bold Z=p\bold Z$. In the following $\nu:K^*@>\quad>>\bold Z$ will be the associated discrete valuation.

\vskip.4cm

\proclaim{\twolverm 8.1 Lemma}\ita Let $\Cal T$ be a diameter four tree of type $(a_1,\dots,a_n)$ such that $\bold Q(\Cal T)^\bullet\subset K$. Suppose that $p$ is $a_i$-regular. Then, there exists a unique normalized model $\beta(X)\in\Cal O_\frak p[X]$ associated to $\Cal T$,  and its reduction modulo $\frak p$ satisfies $\overline\beta(X)=(1-X)^{a_i}\in k(\frak p)[X]$.  \endproclaim

\vskip.4cm

\demo{\twolverm Proof}\es Without any loss of generality, we can assume $i=n$. We will proceed exactly as in the proof of proposition 7.3. First of all, there exists a normalized model defined over $\Cal O_\frak p$. Indeed, if $\beta_0(X)=\prod_{i=1}^n(1-y_iX)^{a_i}$ is any normalized model, let $i\in\{1,\dots,n\}$ such that $\nu(y_i)=\te{Min}\{\nu(y_1),\dots,\nu(y_n)\}$. Then, the polynomial $\beta(X)=\beta_0(y_i^{-1}X)=\prod_{i=1}^n(1-x_iX)^{a_i}$ is a normalized model associated to $\Cal T$, and by construction , we have $x_1,\dots,x_n\in\Cal O_\frak p$, so that $\beta(X)\in\Cal O_\frak p[X]$. Let $N_0=\{i\in\{1,\dots,n-1\}\,\,|\,\,\overline x_i\neq 0\}$.  We can define an equivalence relation on $N_0$ by setting $i\sim j$ if and only if $\overline x_i=\overline x_j$. The quotient $N_1=N_0/\sim$ has cardinality $k\leq n-1$. 
 For any $\tau\in N_1$, set $b_\tau=\sum_{i\in\tau} a_i$ and $x_\tau=\overline x_i$, with $i\in\tau$. Then, $x_\tau\neq x_\sigma$ for any $\tau\neq\sigma$ and the equations in condition iii) of proposition 4.1 give $\sum_{\tau\in N_1}b_\tau x_\tau^i=0$ for any $i\in\{1,\dots,n-1\}$. If we restrict to the first $k$ equations, we obtain a  system of $k$ homogeneous linear equations in the $k$ variables $b_\tau$. Since its determinant does not vanish, we deduce that $b_\tau=0$ for any $\tau\in N_1$, which is absurd, since $b_\tau$ divides $d_n(a_1,\dots,a_n)$. We then deduce that $N_0=N_1=\emptyset$.  Since $\beta(X)$ is a normalized model, this last relation is equivalent to  $x_n=1$ and $x_j\in\frak p$ for any $j<n$. In particular, $\overline\beta(X)=(1-X)^{a_n}$, and no other model is defined over $\Cal O_\frak p$. The lemma is thus proved\qed
\enddemo

\vskip.4cm

The polynomial $\beta(X)$ of the previous lemma will be called the {\bbol canonical model} associated to $\Cal T$. The next result gives the $\frak p$-adic distance between the roots of the canonical model:

\vskip.4cm

\proclaim{\twolverm 8.2 Lemma}\ita The hypothesis being as in the previous lemma, set $\beta(X)=\prod_{j=1}^n(1-x_jX)^{a_j}\in\Cal O_\frak p[X]$. Then $\nu(x_j-x_i)=0$  and $(n-1)\nu(x_j-x_k)=(n-1)\nu(x_j)=\nu(a_i)$ for any $j,k\in\{1,\dots,\widehat i,\dots,n\}$.\endproclaim
 
\vskip.4cm

\demo{\twolverm Proof}\es As in the proof of lemma 8.1, we can assume $i=n$. We know that $x_n=1$ and that $x_1,\dots,x_{n-1}\in\frak p$. In particular, $\nu(x_j-x_n)=0$ for any $j<n$. For $i=n$, the relations in condition iv) of proposition 4.1 give $a_n\prod_{j<n}(x_j-x_n)=(a_n+\dots+a_n)\prod_{j<n}x_j$. Since $p$ does not divide $a_1+\dots+a_n$, we obtain $\nu(a_n)=\sum_{j<n}\nu(x_j)\geq(n-1)\nu$, where $\nu=\te{Min}\{\nu(x_1),\dots,\nu(x_{n-1})\}$.  The equality holds if and only if $\nu=\nu(x_1)=\dots=\nu(x_{n-1})$. Suppose that $(n-1)\nu<\nu(a_n)$ and set $y_j=\pi^{-\nu}x_j\in\Cal O_\frak p$ for any $j<n$ and $y_n=x_n=1$, where $\pi$ is a uniformizer for $\frak p$. Then, the system of equations in condition iii) of proposition 4.1 becomes $a_1y_1^m+\dots+a_{n-1}y_{n-1}^m+\pi^{-m\nu}a_n=0$ for any $m\in\{1,\dots,n-1\}$. Since, by assumption, $(n-1)\nu<\nu(a_n)$, we obtain  $a_1\overline y_1^m+\dots+a_{n-1}\overline y_{n-1}^m=0$ for any $m\leq n-1$. As in the proof of lemma 8.1, set $N_0=\{j\in\{1,\dots,n-1\}\,\,|\,\,\overline y_j\neq 0\}\neq\emptyset$ since, by construction,  there exists $j\in\{1,\dots,n-1\}$ such that $\nu(y_j)=0$.   We can define an equivalence relation on $N_0$ by setting $j\sim k$ if and only if $\overline y_j=\overline y_k$. The quotient $N_1=N_0/\sim$ has cardinality $s\leq n-1$.  For any $\tau\in N_1$, set $b_\tau=\sum_{i\in\tau} a_i$ and $y_\tau=\overline y_i$, with $i\in\tau$. Then, $y_\tau\neq y_\sigma$ for any $\tau\neq\sigma$, and the last  equations  give $\sum_{\tau\in N_1}b_\tau y_\tau^m=0$ for any $m\in\{1,\dots,n-1\}$. If we restrict to the first $s$ equations, we obtain a  system of  homogeneous linear equations in the $s$ variables $b_\tau$. Since its determinant does not vanish, we have $b_\tau=0$ for any $\tau\in N_1$, which is absurd, since $b_\tau$ divides $d_n(a_1,\dots,a_n)$. We then deduce that $(n-1)\nu=\nu(a_n)$. In particular, $\nu(y_j)=0$ for any $j\leq n$. Finally, for $j<n$, the relations in condition iv) of proposition 4.1 give $a_j\prod_{k\neq j}(y_k-y_j)=(a_1+\dots+a_n)\prod_{k\neq j}y_k$, so that, in terms of valuations, since $p$ does not divide $a_j$ nor $a_1+\dots+a_n$, we obtain $\nu(y_k-y_j)=0$. This last relation being equivalent to $\nu(x_k-x_j)=\nu(x_j)$, the lemma is proved\qed 

\enddemo

\vskip.4cm

This last result allows us to obtain some precise results concerning the ramification in the field of moduli of a diameter four tree. First of all, we will need to introduce some  definitions: let $a_1\leq\dots\leq a_n$ be positive integers and suppose tat $p$ is a $a_i$-regular prime. We can define an equivalence relation on  $N_i=\{1,\dots,\widehat i,\dots,n\}$ by setting $j\sim k$ if and only if $a_j=a_k$. We then obtain a decomposition of $N_i$ in $r\leq n-1$ disjoint classes $s_1,\dots,s_r$, and $n_1,\dots,n_r$ will denote their cardinalities. Consider the positive integer $n_0$ defined by 
$$n_0\bold Z=\sum_{j=1}^rn_j\bold Z,$$
that is, $n_0$ is the greatest common divisor of  $n_1,\dots,n_r$. Set $a_i=p^hm$, with $h>0$ and $(m,p)=1$. The integer 
$$e_i=e_i(a_1,\dots,a_n,p)={n-1\over (n-1,hn_0)}$$ 
will be called the {\bbol combinatoral ramification index} associated to $a_i$ and $p$.

\vskip.4cm

\proclaim{\twolverm 8.3 Theorem}\ita Let $\Cal T$ be a diameter four tree of type $(a_1,\dots,a_n)$ over $\QQ$ and suppose that  $p$ is a $a_i$-regular prime. Consider a maximal ideal $\frak p$ of the ring of integers of $\bold Q(\Cal T)$ lying above $p\bold Z$ and denote by $e_\frak p$ its global ramification index. Then, the combinatorial ramification index $e_i$ divides  $e_\frak p$. 
\endproclaim

\vskip.4cm

\demo{\twolverm Proof}\es Consider the canonical model $\beta(X)=\prod_{j=1}^n(1-x_jX)^{a_j}$ associated to $\Cal T$ and let  $s_1,\dots,s_r\subset N_i$ be the equivalence classes defined above. For any distinct $j\in\{1,\dots,r\}$, set $t_j=\prod_{u\in s_j}x_u$. One easily checks that $t_j$  belong to $\bold Q(\Cal T)$. Moreover, lemma 8.2 gives $(n-1)\nu(t_j)=n_j\nu(a_i)$ where $\nu:\bold Q(\Cal T)^*@>\quad>>\bold Z$ is the valuation associated to $\frak p$, normalized by $\nu(\pi)=1$ for any uniformized $\pi$ of $\frak p$. Moreover, we have $\nu(a_i)=\nu(p^hm)=he_\frak p$. Let $m_1,\dots,m_r$ be integers satisfying $n_0=\sum_{j=1}^rm_jn_j$.
Setting $t_0=\prod_{j=1}^rt_j^{m_j}$, we clearly have $t_0\in\bold Q(\Cal T)$ and $(n-1)\nu(t_0)=n_0\nu(a_i)=hn_0e_\frak p$. Finally, let $a,b\in\bold Z$ such that $(n-1,hn_0)=a(n-1)+bhn_0$ and set $t=p^at_0^b\in\bold Q(\Cal T)$. We then obtain 
$$(n-1)\nu(t)=(n-1)\nu(p^a)+(n-1)\nu(t_0^b)=(a(n-1)+bhn_0)e_\frak p=(n-1,hn_0)e_\frak p,$$
from which it follows that $e_i\nu(t)=e_\frak p$. Now, since $\nu(t)$ is an integer, we deduce that $e_i$ divides $e_\frak p$, which concludes the proof.\qed 
\enddemo

\vskip.4cm

We will now describe a quite surprising construction, which will allow us to reduce the study of diameter four trees to the case of good reduction (for $a_i$-regular primes). In the previous paragraphs, we investigated the properties of normalized models. We will now introduce another class of standard model. From now on, we will suppose that $p>n$. In this case, we know from proposition 4.1 that  the polynomial  $\beta(X)=\prod_{i=1}^n(1-x_iX)^{a_i}$ is a standard model associated to a diameter four tree $\Cal T$ of type $(a_1,\dots,a_n)$ if and only if the elements $x_1,\dots,x_n$ satisfy the system of algebraic equations $\phi_1=\dots=\phi_{n-1}=0$, where $\phi_m=\sum_{i=1}^na_ix_i^m$. If these elements verify the further relation $\phi_n=\sum_{i=1}^na_ix_i^n=1$, then we will say that $\beta(X)$ is a {\bbol Kummer model} associated to $\Cal T$. Remark that the expression of $\beta'(X)/\beta(X)$ given at the end of the proof of proposition 4.1 implies that $\phi_n=(-1)^{n-1}x_1\dots x_n(a_1+\dots+a_n)$. As in the case of normalized models, the diameter four tree $\Cal T$ possesses $n\over m$ Kummer models, where $m$ is the order of the automorphism group of $\Cal T$. They can be obtained in the following way: starting from any standard (or normalized) model $\beta(X)=\prod_{i=1}^n(1-x_iX)^{a_i}$ associated to $\Cal T$,  consider the polynomial $h(X)=X^n-t(\beta)$, where $t(\beta)=(-1)^{n-1}x_1\dots x_n(a_1+\dots+a_n)\neq 0$. Then, for any root $x$ of $h(X)$, the polynomial $\beta(x^{-1}X)$ is a Kummer model associated to $\Cal T$. Remark that Hilbert '90 implies that we may choose $\beta(X)\in\bold Q(\Cal T)$. In this case, the splitting field $\bold Q(\Cal T)^\circ$ of $h(X)$ (as an extension of $\bold Q(\Cal T)$) only depends on the diameter four tree $\Cal T$ and not on the standard model $\beta(X)$. Moreover, any Kummer model is defined over $\bold Q(\Cal T)^\circ$. In particular, if $\beta_1(X)$ and $\beta_2(X)$ are two Kummer models associated to $\Cal T$, then there exists a $n$-th root of unity $\zeta\in\mu_n$ such that $\beta_2(X)=\beta_1(\zeta X)$. The field $\bold Q(\te{IV}_{a_1,\dots,a_n})^\circ$ is defined as the compositum of the fields $\bold Q(\Cal T)^\circ$, for all diameter four trees $\Cal T$ of type $(a_1,\dots,a_n)$. It is a Galois extension of $\bold Q$. We will now give a lifting result for Kummer models, which is the analogue of proposition 6.1.

\vskip.4cm

\proclaim{\twolverm 8.4 Proposition}\ita Let $(R,\frak p)$ denote a complete discrete valuation ring of field  of fraction $K$ of characteristic $0$ and of residue field $k=R/\frak p$ of characteritic $p>n$. Suppose that $\overline\beta(X)\in k[X]$ is a Kummer  model associated to a diameter four tree of type $(a_1,\dots,a_n)$. Then, it can be uniquely lifted to a  Kummer model  $\beta(X)\in R[X]$ for a diameter four tree  of the same type. 
\endproclaim

\vskip.4cm

\demo{\twolverm Proof}\es Consider a Galois extension $R'$ of $R$ of field of fraction $K'$ and residue field $k'$, such that $\bold F_p(\Cal T_p)^\bullet\subset k'$.  
 Set $\overline\beta(X)=\prod_{i=1}^n(1-\overline x_iX)^{a_i}$, with $\overline x_1,\dots,\overline x_n\in (k')^*$ pairwise distinct.
Following proposition 4.1, and since $p>n$, the  elements $\overline x_1,\dots,\overline x_{n}$ define a solution of  the system of $n$ algebraic equations $\phi_1=\dots=\phi_{n-1}=0$ and $\phi_n=1$, where $\phi_m=\sum_{i=1}^n\prod_{i=1}^na_iX_i^{m}$. The multidimensional Hensel's lemma then implies that if the determinant of jacobian matrix $J_\phi=(\partial\phi_i/\partial X_j)$ is a unit in $k'$ for $X_1=\overline x_1,\dots, X_n=\overline x_n$,  then this solution  can be uniquely lifted to a solution $x_1,\dots,x_{n-1},x_n$ with $x_i\in R'$ for any $i\in\{1,\dots,n-1\}$. Now, we have $\partial\phi_i/\partial X_j=ia_jX_j^{i-1}$, so that 
$$\te{det}(J_\phi)=\te{det}\left(\matrix a_1 & \dots  & a_n\\ \vdots & \, & \vdots\\ na_1X_1^{n-1}&\dots & na_nX_n^{n-1}\endmatrix\right)=$$
$$=n!a_1\dots a_n\te{det}\left(\matrix 1 & \dots  & 1\\ \vdots & \, & \vdots\\ X_1^{n-1}&\dots & X_n^{n-1}\endmatrix\right)=n!a_1\dots a_n\prod_{1\leq i<j\leq n}(X_i-X_j),$$
which is different from zero for $X_1=\overline x_1,\dots,X_n=\overline x_n$, since $\overline x_i\neq\overline x_j$ for any $i\neq j$ (by hypothesis) and  $p$ does not divide  $a_1\dots a_n$ (the ramification being tame) nor $n!$, because $p>n$. The corresponding Kummer model in characteristic zero is given by $\beta(X)=\prod_{i=1}^n(1-x_iX)^{a_i}$. Let now $\sigma\in\te{Gal}(K'/K)=\te{Gal}(k'/k)$. The polynomial $^\sigma\beta(X)$ obtained by acting on the coefficients of $\beta(X)$ reduce to $\overline\beta(X)$ modulo $\frak p$ (since $\overline\beta(X)\in k[X]$). By uniqueness of the lifting, we then deduce that $^\sigma\beta(X)=\beta(X)$, so that $\beta(X)\in R[X]$.  
\qed  
\enddemo

\vskip.4cm 

\es We will now come back to the theme of this section. Let $\Cal T$ be a diameter four tree of type $(a_1,\dots,a_n)$ and suppose that $p$ is a $a_i$-regular prime. Consider a number field $K$ containing $\bold Q(\Cal T)^\bullet$ and the splitting field of the polynomial $h(X)=X^{n-1}+a_i$. As usual, $\frak p$ will denote a maximal ideal of $\Cal O_K$ lying above $p\bold Z$ and $\nu:K^*@>\quad>>\bold Z$ will be the associated valuation. Let $\beta(X)=\prod_{i=1}^n(1-x_iX)^{a_i}\in\Cal O_\frak p[X]$ be the canonical model associated to $\Cal T$. Fix a root $x$ of $h(X)$ and set $y_j=x^{-1}x_j$ for any $j\neq i$, and $y_i=x_i=1$. Lemma 8.2 impies that $y_1,\dots,\widehat y_i,\dots,y_n$ are elements of $\Cal O_\frak p$ and that they specialise to pairwise distinct elements of $k(\frak p)^*$. We are assuming $p>n$, so that, following proposition 4.1, the elements $x_1,\dots,x_n$ satisfy the system of algebraic equations $\phi_1=\dots=\phi_{n-1}=0$, and this condition is equivalent to the fact that $\beta(X)$ is a normalized model. If we replace $x_1,\dots,x_n$ by $y_1,\dots,y_n$, we then obtain the system $\varphi_1=\dots\varphi_{n-1}=0$, where
$$\varphi_m=x^{-m}\phi_m=\sum_{j\neq i}a_jy_j^m+x^{-m}a_i$$
for any $m\in\{1,\dots,n-1\}$. Moreover, we have $(n-1)\nu(x)=\nu(a_i)$, so that $x^{-m}a_i\in\Cal O_\frak p$. The proof of the following result being  the same as in proposition 8.4, it will be omitted.

\vskip.4cm 

\proclaim{\twolverm 8.5 Proposition}\ita Let $(R,\frak p)$ denote a complete discrete valuation ring of field  of fraction $K$ of characteristic $0$ and of residue field $k=R/\frak p$ of characteritic $p>n$. Suppose that $\overline y_1,\dots,\widehat{\overline y_i},\dots,\overline y_n$ are pairwise distinct elements of $k^*$ satisfying the system of algebraic equations $\varphi_1=\dots=\varphi_{n-1}=0$. Then they can be uniquely lifted to a solution $y_1,\dots,\widehat y_i,\dots,y_n\in R$ of the same system in characteristic zero. 
\endproclaim

\vskip.4cm

\es In order to state the next result, we just need another definition: the fixed root $x$ of $h(X)=X^{n-1}+a_i$ defines a cocycle $\zeta_x\in\te{C}^1(\QQ,\mu_{n-1})$, explicitly given by $\zeta_x(\sigma)=\sigma(x)x^{-1}$ for any $\sigma\in G_\bold Q$. We can twist the action of $G_\bold Q$  by setting $\widehat{\sigma}(t)=\zeta_x(\sigma)\sigma(t)$, for any $t\in\QQ$.

\vskip.4cm

\proclaim{\twolverm 8.6 Theorem}\ita Let $a_1\leq\dots\leq a_n$ be positive integers and suppose that $p>n$ is a $a_i$-regular prime. Consider a number field $K$ containing the fields $\bold Q(\te{IV}_{a_1,\dots,a_n})^\bullet$, $\bold Q(\te{IV}_{a_1,\dots,\widehat a_i,\dots,a_n})^\circ$ and the splitting field of the polynomial $h(X)=X^{n-1}+a_i$. Fix a root $x$ of $h(X)$ and a prime $\frak p$ of $\Cal O_K$ lying over $p\bold Z$. Then, there is a bijection $\Phi=\Phi_{\frak p,x}$ between the set of canonical models associated to diameter four trees of type $(a_1,\dots,a_n)$ and the set of Kummer models associated to diameter four trees of type $(a_1,\dots,\widehat a_i,\dots,a_n)$. Moreover, by construction, all the elements of $\te{IV}_{a_1,\dots,\widehat a_i,\dots,a_n}(\QQ)$ have good reduction at $\frak p$ and, for any $\sigma\in D_\frak p$, we have $\Phi({^\sigma}\beta(X))={^{\widehat\sigma}}\Phi(\beta(X))$.

\endproclaim

\vskip.4cm

\demo{\twolverm Proof}\es Let $\beta(X)=\prod_{j=1}^n(1-x_jX)^{a_j}$ be the canonical model  associated to a diameter four tree of type $(a_1,\dots,a_n)$. As before, if we set $y_j=x^{-1}x_j$ for $j\neq i$ and $y_i=x_i=1$, then $y_1,\dots,y_n$ define a solution of the system of aleraic equations $\varphi_1=\dots=\varphi_{n-1}=0$, where $\varphi_m=\sum_{j\neq i}a_jy_j^m+x^{-m}a_i$. Since $(n-1)\nu(x)=\nu(a_i)$, in $k(\frak p)$ this system can be written as $\sum_{j\neq i}a_j\overline y_j^m=0$ for $m<n-1$ and $\sum_{j\neq i}a_j\overline y_j^{n-1}=1$. But these are just the equations defining Kummer models associated to diameter four trees of type $(a_1,\dots,\widehat a_i,\dots,a_n)$. In particular, $\beta(X)=\prod_{j\neq i}(1-\overline y_j)^{a_j}$ is a Kummer model over $\Fp$ and, following proposition 8.4, it can be uniquely lifted in characteristic zero.

Conversely, if $\beta(X)=\prod_{j\neq i}(1-y_jX)^{a_j}\in K[X]$ is a Kummer model, since $p$ does not divde $d_i(a_1,\dots,a_n)=d(a_1,\dots,\widehat a_i,\dots,a_n)$, we see that $\beta(X)$ (or the corresponding diameter four tree) ha good reduction at $\frak p$, and thus, the elements $\overline y_1,\dots,\overline y_{i-1},\overline y_{i+1},\dots,\overline y_n\in k(\frak p)^*$ are pairwise distinct and define a solution of the system $\phi_1=\dots=\phi_{n-2}=0$ and $\phi_{n-1}=1$. As we noticed,  in $k(\frak p)$ these equations coincide with $\varphi_1=\dots=\varphi_{n-1}=0$. By proposition 8.5, we can uniquely lift the solution $\overline y_1,\dots,\overline y_i=1,\dots,\overline y_n$ of this last system in characteristic zero, and thus  obtain a canonical model associated to a diameter four tree of type $(a_1,\dots,a_n)$. Concerning the Galois action, we just have to remark that, for any $j\neq i$, we have $x^{-1}\sigma(x_j)=x^{-1}\sigma(xy_j)=\sigma(x)x^{-1}\sigma(y_j)=\zeta_x(\sigma)\sigma(y_j)=\widehat\sigma(y_j)$.\qed
\enddemo    

\vskip.4cm 

\proclaim{\twolverm 8.7 Corollary}\ita Let $\Cal T$ be a diameter four tree of type $(a_1,\dots,a_n)$ over $\QQ$ and suppose that $p>n$ is a $a_i$-regular prime. Set $a_i=p^hm$, with $(p,m)=1$ and consider a maximal ideal $\frak p$ of $\Cal O_{\bold Q(\Cal T)}$ lying over $p\bold Z$. If $e_\frak p$ denotes the global index of ramification of $\frak p$, then we have $e_\frak p\leq {n-1\over(n-1,h)}$.

\endproclaim

\vskip.4cm

\demo{\twolverm Proof}\es Fix a root $x$ of $h(X)=X^{n-1}+a_i$ and consider a Galois extension  $K$ of $\bold Q$ satisfying the conditions of theorem 8.6. Let $\frak q$ be a prime of $\Cal O_K$ lying over $\frak p$ and set $\frak r=\frak q\cap\bold Q(x)$. We can easily prove that the extension $\bold Q(x)/\bold Q$ is totally ramified at $\frak r$, and that its ramification index is  $n-1\over(n-1,h)$. Let $\sigma\in\te{Gal}(K/\bold Q(x))$ lying in the inertia subgroup $I_{\frak q/\frak r}$. If $\beta(X)$ is the canonical model associated to $\Cal T$, we have $\Phi({^\sigma}\beta(X))={^{\widehat\sigma}}\Phi(\beta(X))={^\sigma}\Phi(\beta(X))$, since $\sigma$ acts trivially on $x$. Moreover, we know from proposition 7.3 that the inertia at $\frak q$ acts trivially on $\Phi(\beta(X))$ (in fact, we only know that the inertia acts trivially on the set of normalized models, but, in the case of good reduction, we can easily obtain the same result for Kummer models). Thus, we have ${^\sigma}\beta(X)=\beta(X)$, so that $\beta(X)$ is defined over $K^{I_{\frak q/\frak r}}$, which is unramified above $\frak r$, and the corollary follows from this last property.\qed

\enddemo

\vskip.4cm

\example{\twolverm 8.8 examples}

\vskip.2cm

\es Let's come back to  diameter four trees of type $(a_1,\dots,a_n)=(1,\dots,1,a,b)$ introduced in example 4.3.2 (with $a<b$). We already know that this valency class has cardinality $n-1$. We will only treat the case of $b$-regular primes (the $a$-regular case  can be studied in the same way). One easily checks that the integers $d_{n}(a_1,\dots,a_n)=d(a_1,\dots,a_{n-1})=d(1,\dots,a)$ and $c(n,a)={(n-2)!(a+n-2)!\over (a-1)!}$ have the same set of prime divisors. In particular, $p$ is $b$-regular if and only if it divides $b$ and  does not divide $c(n,a)$. This is clearly the case if $p>a+n-2$. The combinatorial ramification index associated to $b$ and $p$ is $e_n={n-1\over(n-1,h)}$, where $p^h$ is the greatest power of $p$ dividing $b$. If $h=1$, i.e., $p$ is a simple divisor of $b$, then theorem 8.3 implies that $e_n=n-1$ divides the index of ramification of any prime of (the ring of integers of) $\bold Q(\Cal T)$ lying above $p\bold Z$, for any tree $\Cal T\in\te{IV}_{1,\dots,1,a,b}(\QQ)$. Now, we have $n-1=e_\frak p\leq[\bold Q(\Cal T):\bold Q]\leq n-1$, so that $\te{IV}_{1,\dots,1,a,b}(\QQ)$  is actually a Galois orbit, and there exist only one prime $\frak p$ of $\bold Q(\Cal T)$ lying above $p\bold Z$. Moreover, $\frak p$ is totally ramified. As an concrete example,  take $n=5$, $a=2$ and $b=7m$, with $m>0$ not divisible by $7$. Since $c(5,2)=2^4\cdot3^2\cdot5$, we see that $p=7$ is $7m$-regular. In particular, $\te{IV}_{1,1,1,2,7m}(\QQ)$ is always a Galois orbit. On the opposite direction, if $n-1$ divides $h$ and $p>n$, then we cannot affirm that  $\Cal T\in\te{IV}_{1,\dots,1,a,b}(\QQ)$ is a Galois orbit, but theorem 8.3 and corollary 8.7 imply that $p\bold Z$ does not ramify in the field of moduli of the diameter four trees of this type (take for example $n=5$, $a=2$, $b=7^4$ and $p=7$).

\endexample

\vskip1.2cm

\specialhead{\twoplverm 9 Wild ramification above infinity}\endspecialhead

\vskip.4cm

\es In this last section, we will study the specialisation of normalized models associated to diameter four trees of tipe $(a_1,\dots,a_n)$ at the primes dividing the degree $a_1+\dots+a_n$ of the covering. The results and techniques are similar to those of the previous section. First of all, we will say that a prime number $p$ is {\bbol regular at infinity} if it divides $a_1+\dots+a_n$ but does not divides the integer $d_\infty(a_1,\dots,a_n)$ defined by $$d_\infty(a_1,\dots,a_n)={1\over a_1+\dots+a_n}d(a_1,\dots,a_n)=\prod_{S\subsetneq\{1,\dots,n\}}\sum_{j\in S}a_j$$
 As usual, $K$ will denote a number field, $\frak p$ a maximal ideal of $K$ lying over $p\bold Z$ and $\nu:K^*@>\quad>>\bold Z$ will be the associated discrete valuation.

\vskip.4cm

\proclaim{\twolverm 9.1 Lemma}\ita Let $\Cal T$ be a diameter four tree of type $(a_1,\dots,a_n)$ such that $\bold Q(\Cal T)^\bullet\subset K$. Suppose that $p$ is regular at infinity. Then, any  normalized model $\beta(X)$ associated to $\Cal T$ is defined over $\Cal O_\frak p$,  and its reduction modulo $\frak p$ satisfies $\overline\beta(X)=(1-X)^{a_1+\dots+a_n}\in k(\frak p)[X]$.  \endproclaim

\vskip.4cm

\demo{\twolverm Proof}\es  We will proceed exactly as in the proof of propositions 7.3 and 8.1. First of all, there exists a normalized model defined over $\Cal O_\frak p$. Indeed, if $\beta_0(X)=\prod_{i=1}^n(1-y_iX)^{a_i}$ is any normalized model, let $i\in\{1,\dots,n\}$ such that $\nu(y_i)=\te{Min}\{\nu(y_1),\dots,\nu(y_n)\}$. Then, the polynomial $\beta(X)=\beta_0(y_i^{-1}X)=\prod_{i=1}^n(1-x_iX)^{a_i}$ is a normalized model associated to $\Cal T$, and by construction , we have $x_1,\dots,x_n\in\Cal O_\frak p$, so that $\beta(X)\in\Cal O_\frak p[X]$. Let $N_0=\{i\in\{1,\dots,n\}\,\,|\,\,\overline x_i\neq 0\}$. Since $\beta(X)$ is  normalized, we have $N_0\neq\emptyset$.  We can define an equivalence relation on $N_0$ by setting $i\sim j$ if and only if $\overline x_i=\overline x_j$.The relations in condition iv) of proposition 4.1 imply the quotient $N_1=N_0/\sim$ has cardinality $k\leq n-1$. 
 For any $\tau\in N_1$, set $b_\tau=\sum_{i\in\tau} a_i$ and $x_\tau=\overline x_i$, with $i\in\tau$. Then, $x_\tau\neq x_\sigma$ for any $\tau\neq\sigma$ and the equations in condition iii) of proposition 4.1 give $\sum_{\tau\in N_1}b_\tau x_\tau^i=0$ for any $i\in\{1,\dots,n-1\}$. If we restrict to the first $k$ equations, we obtain a  system of $k$ homogeneous linear equations in the $k$ variables $b_\tau$. Since its determinant does not vanish, we deduce that $b_\tau=0$ for any $\tau\in N_1$, and our assumptions on $p$ imply that the only possibility is $b_\tau=a_1+\dots+a_n$, so that $k=1$. 
 In particular, since $1\in\{x_1,\dots,x_n\}$, we have $\overline\beta(X)=(1-X)^{a_1+\dots+a_n}$, and the lemma is proved.\qed
\enddemo

\vskip.4cm

\proclaim{\twolverm 9.2 Lemma}\ita The hypothesis being as in the previous lemma, let $\beta(X)=\prod_{i=1}^n(1-x_iX)^{a_i}\in\Cal O_\frak p[X]$ be a normalized model associated to $\Cal T$. Then, for any distinct $i,j\in\{1,\dots,n\}$, we have  $(n-1)\nu(x_i-x_j)=\nu(a_1+\dots+a_n)$.\endproclaim
 
\vskip.4cm

\demo{\twolverm Proof}\es We know that $x_1,\dots,x_{n-1}\in\Cal O_\frak p^*$, since they all reduce to $1$ modulo $\frak p$. Consider an element $i\in\{1,\dots,n\}$. Since $p$ does not divide $a_i$,  the relations in condition iv) of proposition 4.1 lead to $\nu(a_1+\dots+a_n)=\sum_{j\neq i}\nu(x_j-x_i)\geq(n-1)\nu$, where $\nu=\te{Min}_{j\neq i}\{\nu(x_j-x_i)\}$. Suppose that $(n-1)\nu<\nu(a_1+\dots+a_n)$ and set $y_j=\pi^{-\nu}(x_j-x_i)\in\Cal O_\frak p$ for any $j\in\{1,\dots,n\}$, where $\pi$ is a uniformizer for $\frak p$. Then,  the system of equations in condition iii) of proposition 4.1 becomes $a_1y_1^m+\dots+a_ny_n^m=(-1)^mx_i^m\pi^{-m\nu}(a_1+\dots+a_n)$ for any $m\in\{1,\dots,n-1\}$. Since, by assumption, $(n-1)\nu<\nu(a_1+\dots+a_n)$, we obtain  $a_1\overline y_1^m+\dots+a_{n-1}\overline y_{n-1}^m=0$ for any $m\leq n-1$. Set $N_0=\{j\in\{1,\dots,n\}\,\,|\,\,\overline y_j\neq 0\}\neq\emptyset$ since, by construction,  there exists $j\in\{1,\dots,n-1\}$ such that $\nu(y_j)=\nu(x_j-x_i)-\nu=0$.   We can define an equivalence relation on $N_0$ by setting $j\sim k$ if and only if $\overline y_j=\overline y_k$. The quotient $N_1=N_0/\sim$ has cardinality $k\leq n-1$, since $y_i=0$.  For any $\tau\in N_1$, set $b_\tau=\sum_{j\in\tau} a_j$ and $y_\tau=\overline y_j$, with $j\in\tau$. Then, $y_\tau\neq y_\sigma$ for any $\tau\neq\sigma$, and the last  equations  give $\sum_{\tau\in N_1}b_\tau y_\tau^m=0$ for any $m\in\{1,\dots,n-1\}$. If we restrict to the first $k$ equations, we obtain a  system of $k$ homogeneous linear equations in the $k$ variables $b_\tau$. Since its determinant does not vanish, we have $b_\tau=0$ for any $\tau\in N_1$, which is absurd, since $b_\tau$ divides $d_\infty(a_1,\dots,a_n)$. We then deduce that $(n-1)\nu=\nu(a_1+\dots+a_n)$, and thus $(n-1)\nu(x_j-x_i)=\nu(a_1+\dots+a_n)$ for any $j\neq i$. This last relations holding for any $i\in\{1,\dots,n\}$, the lemma is proved.\qed 

\enddemo

\vskip.4cm

As in the previous section, this result has many interesting arithmetical applications.  Let $a_1\leq\dots\leq a_n$ be positive integers and suppose tat $p$ is a regular prime at infinity. We can define an equivalence relation on  $N_0=\{1,\dots,n\}$ by setting $i\sim j$ if and only if $a_i=a_j$. We then obtain a decomposition of $N_0$ in $r\leq n$ disjoint classes $s_1,\dots,s_r$, and $n_1,\dots,n_r$ will denote their cardinalities. Consider the positive integer $n_0$ defined by 
$$n_0\bold Z=\sum_{1\leq i<j\leq r}n_in_j\bold Z+\sum_{i=1}^rn_i(n_i-1)\bold Z,$$
that is, $n_0$ is the greatest common divisor of  $n_in_j$ (for $i<j$) and $n_i(n_i-1)$ (for $1\leq i\leq r$). Set $a_1+\dots+a_n=p^hm$, with $h>0$ and $(m,p)=1$. The integer 
$$e_\infty=e_\infty(a_1,\dots,a_n,p)={n-1\over (n-1,hn_0)}$$ 
will be called the {\bbol combinatoral ramification index} at infinity associated to $p$.

\vskip.4cm

\proclaim{\twolverm 9.3 Theorem}\ita Let $\Cal T$ be a diameter four tree of type $(a_1,\dots,a_n)$ over $\QQ$ and suppose that  $p$ is a regular prime at infinity. Consider a maximal ideal $\frak p$ of the ring of integers of $\bold Q(\Cal T)$ lying above $p\bold Z$ and denote by $e_\frak p$ its global ramification index. Then, the combinatorial ramification index $e_\infty$ divides  $e_\frak p$. 
\endproclaim

\vskip.4cm

\demo{\twolverm Proof}\es Consider a normalized  model $\beta(X)=\prod_{i=1}^n(1-x_iX)^{a_i}$ associated to $\Cal T$ and let  $s_1,\dots,s_r\subset N_0$ be the equivalence classes defined above. For any  $i,j\in\{1,\dots,r\}$, with $i<j$, set $t_i=\prod_{u,v\in s_i,\, u\neq v}(x_u-x_v)$ and $t_{i,j}=\prod_{u\in s_i,\, v\in s_j}(x_u-x_v)$. One easily checks that $t_i$ and $t_{i,j}$ belong to $\bold Q(\Cal T)$. Moreover, lemma 8.2 gives 
$$(n-1)\nu(t_i)=n_i(n_i-1)\nu(a_1+\dots+a_n)\quad\te{and}\quad(n-1)\nu(t_{i,j})=n_in_j\nu(a_1+\dots+a_n),$$
where $n_i$ is the cardinality of $s_i$ and $\nu:\bold Q(\Cal T)^*@>\quad>>\bold Z$ is the valuation associated to $\frak p$, normalized by $\nu(\pi)=1$ for any uniformized $\pi$ of $\frak p$. Moreover, we have $\nu(a_1+\dots+a_n)=\nu(p^hm)=he_\frak p$. Let $m_i,m_{i,j}$ ($i,j\in\{1,\dots,r\}$ and $i<j$) be integers satisfying 
$$n_0=\sum_{1\leq i<j\leq r}m_{i,j}n_in_j+\sum_{i=1}^rm_in_i(n_i-1)$$
Setting $t_0=\prod_{1\leq i<j\leq r}t_{i,j}^{m_{i,j}}\prod_{i=1}^rt_i^{m_i}$, we clearly have $t_0\in\bold Q(\Cal T)$ and $(n-1)\nu(t_0)=n_0\nu(a_1+\dots+a_n)=hn_0e_\frak p$. Finally, let $a,b\in\bold Z$ such that $(n-1,hn_0)=a(n-1)+bhn_0$ and set $t=p^at_0^b\in\bold Q(\Cal T)$. We then obtain 
$$(n-1)\nu(t)=(n-1)\nu(p^a)+(n-1)\nu(t_0^b)=(a(n-1)+bhn_0)e_\frak p=(n-1,hn_0)e_\frak p,$$
from which it follows that $e_\infty\nu(t)=e_\frak p$. Now, since $\nu(t)$ is an integer, we deduce that $e_\infty$ divides $e_\frak p$, which concludes the proof.\qed 
\enddemo

\vskip.4cm

We will finally give the analogue of theorem 8.6.  Let $\Cal T$ be a diameter four tree of type $(a_1,\dots,a_n)$ and suppose that the prime number $p$ is regular at infinity. Consider a number field $K$ containing $\bold Q(\Cal T)^\bullet$ and the splitting field of the polynomial $h(X)=X^{n-1}+(-1)^n(a_1+\dots+a_n)$. As usual, $\frak p$ will denote a maximal ideal of $\Cal O_K$ lying above $p\bold Z$ and $\nu:K^*@>\quad>>\bold Z$ will be the associated valuation. Let $\beta(X)=\prod_{i=1}^n(1-x_iX)^{a_i}\in\Cal O_\frak p[X]$ be a normalized  model associated to $\Cal T$. Suppose that $x_i=1$. In this case, we will say that $\beta(X)$ is a {\bbol $\bold a_i$-normalized model}. Denote by $n(a_i)$ the cardinality of the set $\{j\in\{1,\dots,n\}\,\,|\,\,a_j=a_i\}$. One easily shows that the order $m$ of the automorphism group of $\Cal T$ divides $n(a_i)$ and that there exist exactly $n(a_i)\over m$ $a_i$-normalized models associated to $\Cal T$.  Fix a root $x$ of $h(X)$ and set $y_j=x^{-1}(x_j-1)$ for any $i\in\{1,\dots,n\}$. Lemma 9.2 impies that $y_1,\dots,y_n$ are elements of $\Cal O_\frak p$ and that they specialise to pairwise distinct elements of $k(\frak p)^*$. Moreover, by construction, we have $y_i=0$. We are assuming $p>n$, so that, following proposition 4.1, the elements $x_1,\dots,x_n$ satisfy the system of algebraic equations $\phi_1=\dots=\phi_{n-1}=0$, and this condition is equivalent to the fact that $\beta(X)$ is a normalized model. As in the proof of lemma 9.2, if we replace $x_1,\dots,x_n$ by $y_1,\dots,y_n$, we then obtain the system $\chi_1=\dots\chi_{n-1}=0$, where
$$\chi_m=\sum_{j=1}^na_jy_j^m+(-1)^{m-1}x^{-m}(a_1+\dots+a_n)$$
for any $m\in\{1,\dots,n-1\}$. Moreover, we have $(n-1)\nu(x)=\nu(a_1+\dots+a_n)$, so that $x^{-m}(a_1+\dots+a_n)\in\Cal O_\frak p$. The proof of the following result being  the same as in proposition 8.4, it will be omitted.

\vskip.4cm 

\proclaim{\twolverm 9.4 Proposition}\ita Let $(R,\frak p)$ denote a complete discrete valuation ring of field  of fraction $K$ of characteristic $0$ and of residue field $k=R/\frak p$ of characteritic $p>n$. Suppose that $\overline y_1,\dots,\overline y_n$ are pairwise distinct elements of $k^*$ satisfying the system of algebraic equations $\chi_1=\dots=\chi_{n-1}=0$, with $\overline y_i=0$. Then they can be uniquely lifted to a solution $y_1,\dots,y_n\in R$ of the same system in characteristic zero, with $y_i=0$. 
\endproclaim

\vskip.4cm

\es As in the previous section, a fixed root $x$ of $h(X)=X^{n-1}+(-1)^{n}(a_1+\dots+a_n)$ defines a cocycle $\zeta_x\in\te{C}^1(\QQ,\mu_{n-1})$, explicitly given by $\zeta_x(\sigma)=\sigma(x)x^{-1}$ for any $\sigma\in G_\bold Q$. We can twist the action of $G_\bold Q$  by setting $\widehat{\sigma}(t)=\zeta_x(\sigma)\sigma(t)$, for any $t\in\QQ$. 

\vskip.4cm

\proclaim{\twolverm 9.5 Theorem}\ita Let $a_1\leq\dots\leq a_n$ be positive integers and suppose that the prime number $p>n$ is regular at infinity. Fix an element $i\in\{1,\dots,n\}$ and  consider a number field $K$ containing the fields $\bold Q(\te{IV}_{a_1,\dots,a_n})^\bullet$, $\bold Q(\te{IV}_{a_1,\dots,\widehat a_i,\dots,a_n})^\circ$ and the splitting field of the polynomial $h(X)=X^{n-1}+(-1)^n(a_1+\dots+a_n)$. Fix a root $x$ of $h(X)$ and a prime $\frak p$ of $\Cal O_K$ lying over $p\bold Z$. Then, there is a bijection $\Psi=\Psi_{a_i,\frak p,x}$ between the set of $a_i$-normalized models associated to diameter four trees of type $(a_1,\dots,a_n)$ with $x_i=1$ and the set of Kummer models associated to diameter four trees of type $(a_1,\dots,\widehat a_i,\dots,a_n)$. Moreover, by construction, all the elements of $\te{IV}_{a_1,\dots,\widehat a_i,\dots,a_n}(\QQ)$ have good reduction at $\frak p$ and, for any $\sigma\in D_\frak p$, we have $\Psi({^\sigma}\beta(X))={^{\widehat\sigma}}\Psi(\beta(X))$.

\endproclaim

\vskip.4cm

\demo{\twolverm Proof}\es Let $\beta(X)=\prod_{j=1}^n(1-x_jX)^{a_j}$ be the normalized model  associated to a diameter four tree of type $(a_1,\dots,a_n)$ for which $x_i=1$. As before, if we set $y_j=x^{-1}(x_j-1)$, then $y_1,\dots,y_n$ define a solution of the system of aleraic equations $\chi_1=\dots=\chi_{n-1}=0$, where $\chi_m=\sum_{j=1}^na_jy_j^m+(-1)^{m-1}x^{-m}(a_1+\dots+a_n)$ and $y_i=0$. Since $(n-1)\nu(x)=\nu(a_1+\dots+a_n)$, in $k(\frak p)$ this system can be written as $\sum_{j\neq i}a_j\overline y_j^m=0$ for $m<n-1$ and $\sum_{j\neq i}a_j\overline y_j^{n-1}=1$, with $\overline y_i=0$. But these are just the equations defining Kummer models associated to diameter four trees of type $(a_1,\dots,\widehat a_i,\dots,a_n)$. In particular, $\beta(X)=\prod_{j\neq i}(1-\overline y_jX)^{a_j}$ is a Kummer model over $\Fp$ and, following proposition 8.4, it can be uniquely lifted in characteristic zero.

Conversely, if $\beta(X)=\prod_{j\neq i}(1-y_jX)^{a_j}\in K[X]$ is a Kummer model, since $p$ does not divde $d_\infty(a_1,\dots,a_n)$, we see that $\beta(X)$ (or the corresponding diameter four tree) ha good reduction at $\frak p$, and thus, the elements $\overline y_1,\dots,\overline y_{i-1},\overline y_{i+1},\dots,\overline y_n\in k(\frak p)^*$ are pairwise distinct and define a solution of the system $\phi_1=\dots=\phi_{n-2}=0$ and $\phi_{n-1}=1$. As we noticed,  in $k(\frak p)$ these equations coincide with $\chi_1=\dots=\chi_{n-1}=0$. By proposition 8.5, we can uniqely lift the solution $\overline y_1,\dots,\overline y_i=0,\dots,\overline y_n$ of this last system in characteristic zero, and thus  obtain a $a_i$-normalized model associated to a diameter four tree of type $(a_1,\dots,a_n)$. Concerning the Galois action, we just have to remark that, for any $j\neq i$, we have $x^{-1}\sigma(x_j)=x^{-1}\sigma(xy_j)=\sigma(x)x^{-1}\sigma(y_j)=\zeta_x(\sigma)\sigma(y_j)=\widehat\sigma(y_j)$.\qed
\enddemo    

\vskip.4cm 
 
\proclaim{\twolverm 9.6 Corollary}\ita Let $\Cal T$ be a diameter four tree of type $(a_1,\dots,a_n)$ over $\QQ$ and suppose that the prime number $p>n$ is  regular at infinity. Set $a_1+\dots+a_n=p^hm$, with $(p,m)=1$ and consider a maximal ideal $\frak p$ of $\Cal O_{\bold Q(\Cal T)}$ lying over $p\bold Z$. If $e_\frak p$ denotes the global index of ramification of $\frak p$, then we have $e_\frak p\leq {n-1\over(n-1,h)}$.

\endproclaim  

\vskip.4cm

\example{\twolverm 9.7 examples}

\vskip.2cm

\es As in the previous section, let's study diameter four trees of type $(a_1,\dots,a_n)=(1,\dots,1,a,b)$, with $a<b$.  The integers $d_\infty(a_1,\dots,a_n)=d_\infty(1,\dots,1,a,b)={1\over a+b+n-2}d(1,\dots,1,a,b)$ and $u(n,a,b)={(n-2)!(a+n-2)!(b+n-2)!(a+b+n-3)!\over (a-1)!(b-1)!(a+b-1)!}$ have the same set of prime divisors. In particular, $p$ is regular at infinity if and only if it divides $a+b+n-2$ and  does not divide $u(n,a,b)$.  In this case, the combinatorial ramification index at infinity associated to $p$ is $e_\infty={n-1\over(n-1,h)}$, where $p^h$ is the greatest power of $p$ dividing $a+b+n-2$. If $h=1$,  then theorem 9.3 implies that $e_\infty=n-1$ divides the index of ramification of any prime of (the ring of integers of) $\bold Q(\Cal T)$ lying above $p\bold Z$, for any tree $\Cal T\in\te{IV}_{1,\dots,1,a,b}(\QQ)$. Now, we have $n-1=e_\frak p\leq[\bold Q(\Cal T):\bold Q]\leq n-1$, so that $\te{IV}_{1,\dots,1,a,b}(\QQ)$  is actually a Galois orbit, and there exist only one prime $\frak p$ of $\bold Q(\Cal T)$ lying above $p\bold Z$. Moreover, $\frak p$ is totally ramified. As an concrete example,  take $n=5$, $a=2$ and $b=77$. Since $u(5,2,77)=2^{11}\cdot3^6\cdot5^5\cdot19\cdot37^2\cdot73$, we see that the primes $7$ and $11$ are both  $77$-regular. In particular, $\te{IV}_{1,1,1,2,77}(\QQ)$ a Galois orbit. On the opposite direction, if $n-1$ divides $h$ and $p>n$, then we cannot affirm that  $\Cal T\in\te{IV}_{1,\dots,1,a,b}(\QQ)$ is a Galois orbit, but theorem 9.3 and corollary 9.6 imply that $p\bold Z$ does not ramify in the field of moduli of the diameter four trees of this type (take for example $n=5$, $a=2$, $b=2396=2^2\cdot599$ and $p=7$).

\endexample

\vskip1.2cm

\specialhead{\twoplverm References}\endspecialhead

\vskip.4cm

\es 

\noindent [G] Grothendieck, A. {\ita Esquisse d'un Programme}, in Geometric Galois
Actions 1, London 

\ Math. Soc.  Lecture Note Ser.  242, Cambridge University Press, Cambridge, 1997, 

\vskip.2cm

\noindent [K] Komatsu, T, {\ita Geometric balance of cuspidal points realizing dessins d'enfants on 

\ the Riemann sphere} Math Ann 320 (2001) 3, 417-429.

\vskip.2cm

\noindent [S] Schneps, L., {\ita The Grothendieck theory of dessins d'enfants},
London Math. Soc. 

\ Lecture Note Ser.  200, Cambridge University Press, Cambridge, 1994.

\vskip.2cm

\noindent [Sh] Shabat, G. B. {\ita On the classification of plane trees by their Galois orbits}, in  {\ita The 

\ Grothendieck theory of dessins d'enfants} (Luminy, 1993),
169--177, London Math. 

\ Soc. Lecture Note Ser., 200, Cambridge Univ. Press, Cambridge, 1994. 

\vskip.2cm

\noindent [SV] Shabat, G. B.; Voevodsky, V. A. {\ita Drawing curves over number fields}, in {\ita The 

\  Grothendieck Festschrift}, Vol. III, 199--227, Progr. Math., 88, Birkh\"auser Boston, 

\ Boston, MA, 1990.

\vskip.2cm

\noindent [Z1] Zapponi, L., {\ita Dessins d'enfants et action galoisienne}, Th\`ese de doctorat, Universit\'e 

\ de Franche-Comt\'e, Besan\c con, 1998.

\vskip.2cm 

\noindent [Z2] Zapponi, L., {\ita Fleurs, arbres et cellules: un invariant galoisien pour une famille 

\ d'arbres},  Compositio Math. 122 (2000), no. 2, 113--133.
\end